\documentclass[aos,MSNbibl,nameyear,seceqn,dvips]{arximspdf}
\usepackage{graphicx}
\doi{10.1214/12-AOS999} 
\volume{40}
\issue{4}
\pubyear{2012}
\firstpage{2014}
\lastpage{2042}

\makeatletter

\newcommand{\rright}{\right}
\newcommand{\lleft}{\left}
\newcommand{\rrVert}{\Vert}
\newcommand{\llVert}{\Vert}
\newtheorem{theorem}{Theorem}[section]

\newtheorem{lemma}{Lemma}

\newproclaim{definition}{Definition}
\newproclaim{remark}{Remark}
\newproclaim{assumption}{Assumption}

\def\bx{\mathbf{x}}
\def\bu{\mathbf{u}}
\def\bv{\mathbf{v}}
\def\bE{\mathbb{E}}
\def\bP{\mathbb{P}}
\def\bR{\mathbb{R}}
\def\Bcal{\mathcal{B}}
\def\Ccal{\mathcal{C}}
\def\Fcal{\mathcal{F}}
\def\Ical{\mathcal{I}}
\def\Ncal{\mathcal{N}}
\def\Scal{\mathcal{S}}
\def\card{\operatorname{card}}
\def\ind{\mathbb{I}}

\makeatother

\begin{document}
\begin{frontmatter}

\title{Adaptive covariance matrix estimation through block thresholding}
\runtitle{Adaptive covariance matrix estimation}

\begin{aug}
\author[a]{\fnms{T. Tony} \snm{Cai}\thanksref{t1}\ead[label=e1]{tcai@wharton.upenn.edu}}
\and
\author[b]{\fnms{Ming} \snm{Yuan}\corref{}\thanksref{t2}\ead[label=e2]{myuan@isye.gatech.edu}}
\thankstext{t1}{Supported in part by NSF FRG Grant DMS-08-54973.}
\thankstext{t2}{Supported in part by NSF Career Award DMS-08-46234.}
\runauthor{T. T. Cai and M. Yuan}
\affiliation{University of Pennsylvania and Georgia Institute of Technology}
\address[a]{Department of Statistics\\
The Wharton School\\
University of Pennsylvania\\
Philadelphia, Pennsylvania 19104\\
USA\\
\printead{e1}}

\address[b]{School of Industrial\\
\quad and Systems Engineering\\
Georgia Institute of Technology\\
Atlanta, Georgia 30332\\
USA\\
\printead{e2}}
\end{aug}

\received{\smonth{9} \syear{2011}}
\revised{\smonth{3} \syear{2012}}

%
\begin{abstract}
Estimation of large covariance matrices has drawn considerable
recent attention, and the theoretical focus so far has mainly
been on developing a minimax theory over a fixed parameter space.
In this paper, we consider adaptive covariance matrix estimation
where the goal is to construct a single procedure which is minimax
rate optimal simultaneously over each parameter space in a large collection.
A fully data-driven block thresholding estimator is proposed. The
estimator is
constructed by carefully dividing the sample covariance matrix into blocks
and then simultaneously estimating the entries in a block by thresholding.
The estimator is shown to be optimally rate adaptive over a wide range of
bandable covariance matrices. A~simulation study is carried out and shows
that the block thresholding estimator performs well numerically. Some of
the technical tools developed in this paper can also be of independent interest.
\end{abstract}

%
\begin{keyword}[class=AMS]
\kwd[Primary ]{62H12}
\kwd[; secondary ]{62F12}
\kwd{62G09}
\end{keyword}
\begin{keyword}
\kwd{Adaptive estimation}
\kwd{block thresholding}
\kwd{covariance matrix}
\kwd{Frobenius norm}
\kwd{minimax estimation}
\kwd{optimal rate of convergence}
\kwd{spectral norm}
\end{keyword}

\end{frontmatter}

\section{Introduction}
\label{secint}

Covariance matrix estimation is of fundamental importance in
multivariate analysis. Driven by a wide range of applications in
science and engineering, the high-dimensional setting, where the
dimension $p$ can be much larger than the sample size $n$, is of
particular current interest. In such a setting, conventional methods
and results based on fixed $p$ and large $n$ are no longer applicable,
and in particular, the commonly used sample covariance matrix and
normal maximum likelihood estimate perform poorly.

A number of regularization methods, including banding, tapering,
thresholding and $\ell_1$ minimization, have been developed in recent
years for estimating a large covariance matrix or its inverse. See, for
example, \citet{lw04},\vadjust{\goodbreak} \citet{hlpl06}, \citet{yl07},
\citet{bgd08}, Bickel and Levina (\citeyear{bl08a,bl08b}), \citet{e08},
\citet{ffl08},
Friedman, Hastie and Tibshirani (\citeyear{fht08}), \citet
{rzy08}, \citet{rblz08}, \citet{lf09}, \citet{rlz09},
\citet{czz10}, \citet{y10}, \citet{cl11} and
\citet{cll11}, among many others.

Let $X^{(1)},\ldots, X^{(n)}$ be $n$ independent copies of a $p$
dimensional Gaussian random vector $X=(X_1,\ldots,X_p)^\mathsf{T}\sim
N(\mu
,\Sigma)$. The goal is to estimate the covariance matrix $\Sigma$ and
its inverse $\Sigma^{-1}$ based on the sample $ \{ X^{(i)}\dvtx
i=1,\ldots,n \} $. It is now well known that the usual sample
covariance matrix
\[
\bar{\Sigma}=\frac{1}{n-1}\sum_{i=1}^n
\bigl(X^{(i)}-\bar{X} \bigr) \bigl(X^{(i)}-\bar{X}
\bigr)^\mathsf{T},
\]
where $\bar{X}=\frac{1}{n}\sum_{i=1}^n X^{(i)}$, is not a consistent
estimator of the covariance matrix~$\Sigma$ when $p\gg n$, and
structural assumptions are required in order to estimate $\Sigma$ consistently.

\begin{figure}
\centering
\begin{tabular}{@{}cc@{}}

\includegraphics{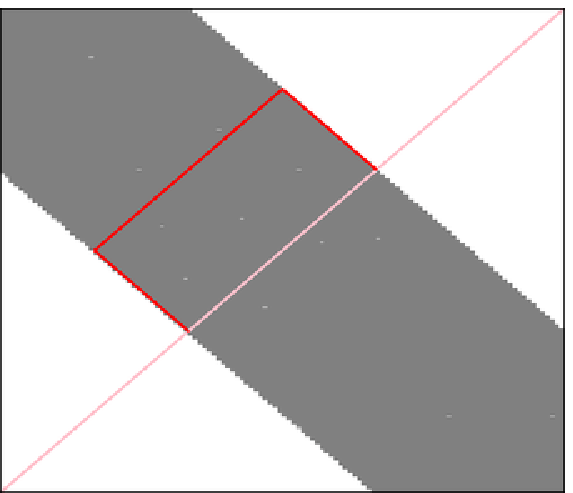}
 & \includegraphics{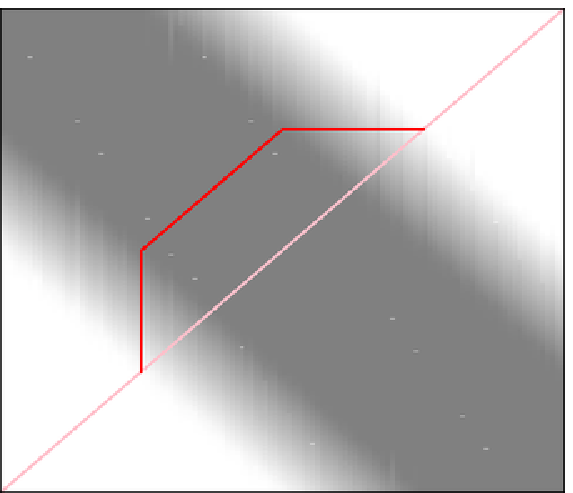}\\
\footnotesize{(a) Weighting matrix for banding} & \footnotesize{(b) Weighting matrix for tapering}\\
\end{tabular}
\caption{Both banding and tapering estimators can be expressed as the
Schur product of the sample covariance matrix and a weighting matrix.
Subfigures of \textup{(a)} and \textup{(b)}
illustrate the weighting matrix for both estimators.}\label{figbandtaper}
\end{figure}

One of the most commonly considered classes of covariance matrices is
the ``bandable'' matrices, where the entries of the matrix decay as they
move away from the diagonal. More specifically, consider the following
class of covariance matrices introduced in \citet{bl08a}:
%
\begin{eqnarray}
\label{paraspace} \Ccal_\alpha&=&\Ccal_\alpha(M_0,M):=
\biggl\{\Sigma\dvtx \max_j\sum_i\bigl\{ |
\sigma_{ij}|\dvtx |i-j|\ge k\bigr\}\le Mk^{-\alpha}\ \forall k,
\nonumber
\\[-8pt]
\\[-8pt]
\nonumber
&&\hspace*{80pt}\mbox{and } 0< M_0^{-1} \le
\lambda_{\min}( \Sigma), \lambda_{\max
}(\Sigma)\le
M_0 \biggr\}.
\end{eqnarray}
Such a family of covariance matrices naturally arises in a number of
settings, including temporal or spatial data analysis. See \citet
{bl08a} for further discussions.
Several regularization methods have been introduced for estimating a
bandable covariance matrix $\Sigma\in\Ccal_\alpha$.
\citet{bl08a} suggested banding the sample covariance
matrix $\bar{\Sigma}$ and estimating $\Sigma$ by $\bar{\Sigma
}\circ
B_k$ where $B_k$ is a banding matrix
\[
B_k= \bigl({\mathbb I} \bigl(|i-j|\le k \bigr) \bigr)_{1\le i,j\le p}
\]
and $\circ$ represents the Schur product, that is, $(A\circ
B)_{ij}=A_{ij}B_{ij}$ for two matrices of the same dimensions.
See Figure~\ref{figbandtaper}(a) for an illustration.
\citet{bl08a} proposed to choose $k\asymp(n/\log
p)^{1/(2(\alpha+1))}$ and showed that the resulting banding estimator
attains the rate of convergence
%
\begin{equation}
\label{BLrate} \|\bar{\Sigma}\circ B_k-\Sigma\|=O_p
\biggl( \biggl(\frac{\log p}{n} \biggr)^{{\alpha}/{(2\alpha+2)}} \biggr)
\end{equation}
uniformly over $\Ccal_\alpha$, where $\|\cdot\|$ stands for the
spectral norm.
This result indicates that even when $p\gg n$, it is still possible to
consistently estimate $\Sigma\in\Ccal_\alpha$, so long as $\log
p=o(n)$.\vadjust{\goodbreak}

\citet{czz10} established the minimax rate of convergence
for estimation over $\Ccal_\alpha$ and introduced a tapering estimator
$\bar{\Sigma}\circ T_k$ where the tapering matrix $T_k$ is given by
\[
T_k= \biggl(\frac{2}{k} \bigl\{\bigl(k-|i-j|\bigr)_+-\bigl(k/2-|i-j|\bigr)_+
\bigr\} \biggr)_{1\le i,j\le p},
\]
with $(x)_+=\max(x,0)$. See Figure~\ref{figbandtaper}(b) for an illustration. It was shown that the tapering estimator $\bar
{\Sigma}\circ T_k$ with $k\asymp n^{1/(2\alpha+1)}$ achieves the rate
of convergence
%
\begin{equation}
\label{optrate} \|\bar{\Sigma}\circ T_k-\Sigma\|=O_p
\biggl(n^{-{\alpha}/{(2\alpha
+1)}}+ \biggl(\frac{\log p}{n} \biggr)^{1/2} \biggr)
\end{equation}
uniformly over $\Ccal_\alpha$, which is always faster than the rate in
(\ref{BLrate}).
This implies that the rate of convergence given in~(\ref{BLrate}) for
the banding estimator with $k\asymp(n/\log p)^{1/(2(\alpha+1))}$ is in
fact sub-optimal.
Furthermore, a lower bound argument was given in \citet{czz10}
which showed that the rate of convergence given in (\ref
{optrate}) is indeed optimal for estimating the covariance matrices
over~$\Ccal_\alpha$.

The minimax rate of convergence in~(\ref{optrate}) provides an
important benchmark for the evaluation of the performance of covariance
matrix estimators. It is, however, evident from its construction that
the rate optimal tapering estimator constructed in \citet{czz10}
requires explicit knowledge of the decay rate $\alpha$ which is
typically unknown in practice. It is also clear that a tapering
estimator designed for a parameter space with a given decay rate
$\alpha
$ performs poorly over another parameter space with a different decay
rate. The tapering estimator mentioned above is thus not very practical.\vadjust{\goodbreak}

This naturally leads to the important question of adaptive estimation:
Is it possible to construct a single estimator, not depending on the
decay rate $\alpha$, that achieves the optimal rate of convergence
simultaneously over a wide range of the parameter spaces $\Ccal_\alpha$?
We shall show in this paper that the answer is affirmative.
A fully data-driven adaptive estimator $\hat\Sigma$ is constructed and
is shown to be simultaneously rate optimal over the collection of the
parameter spaces $\Ccal_\alpha$ for all $\alpha>0$. That is,
\[
\sup_{\Sigma\in\Ccal_\alpha}\bE\|\hat{\Sigma}-\Sigma\|^2\asymp \min \biggl
\{n^{-{2\alpha}/{(2\alpha+1)}}+\frac{\log p}{n}, \frac{p}{n} \biggr\} \qquad\mbox{for all $
\alpha>0$}.
\]
In many applications, the inverse covariance matrix is of significant
interest. We introduce a slightly modified version of $\hat\Sigma^{-1}$
and show that it adaptively attains the optimal rate of
convergence for estimating $\Sigma^{-1}$.


The adaptive covariance matrix estimator achieves its adaptivity
through block thresholding of the sample covariance matrix $\bar\Sigma$.
The idea of adaptive estimation via block thresholding can be traced
back to nonparametric function estimation using Fourier or wavelet
series. See, for example, \citet{e85} and \citet{c99}. However,
the application of block thresholding to covariance matrix estimation
poses new challenges. One of the main difficulties in dealing with
covariance matrix estimation, as opposed to function estimation or
sequence estimation problems, is the fact that the spectral norm is not
separable in its entries. Another practical challenge is due to the
fact that the covariance matrix is ``two-directional'' where one
direction is along the rows and another along the columns. The blocks
of different sizes need to be carefully constructed so that they fit
well in the sample covariance matrix and the risk can be assessed based
on their joint effects rather than their individual contributions.
There are two main steps in the construction of the adaptive covariance
matrix estimator. The first step is the construction of the blocks.
Once the blocks are constructed, the second step is to estimate the
entries of the covariance matrix $\Sigma$ in groups and make
simultaneous decisions on all the entries within a block. This is done
by thresholding the sample covariance matrix block by block. The
threshold level is determined by the location, block size and
corresponding spectral norms. The detailed construction is given in
Section~\ref{secmeth}.

We shall show that the proposed block thresholding estimator $\hat
\Sigma$ is simultaneously rate-optimal over every $\Ccal_\alpha$ for
all $\alpha>0$.
The theoretical analysis of the estimator $\hat\Sigma$ requires some
new technical tools that can be of independent interest. One is a
concentration inequality which shows that although the sample
covariance matrix $\bar{\Sigma}$ is not a reliable estimator of
$\Sigma
$, its submatrices could still be a good estimate of the corresponding
submatrices of $\Sigma$.
Another useful tool is a so-called norm compression inequality which
reduces the analysis on the whole matrix to a matrix of much smaller
dimensions, whose entries are the spectral norms of the blocks.


In addition to the analysis of the theoretical properties of the
proposed adaptive block thresholding estimator, a simulation study is
carried out to investigate the finite sample performance of the
estimator. The simulations show that the proposed estimator enjoys good
numerical performance when compared with nonadaptive estimators such as
the banding and tapering estimators.

Besides bandable matrices considered in the present paper, estimating
sparse covariance matrices and sparse precision matrices has also been
actively studied in the recent literature. \citet{bl08b}
proposed a thresholding estimator for sparse covariance matrices and
obtained the rate of convergence. \citet{cz11} developed a new
general minimax lower bound technique and established the minimax rate
of convergence for estimating sparse covariance matrices under the
spectral norm and other matrix operator norms. \citet{cl11}
introduced an adaptive thresholding procedure for estimating sparse
covariance matrices that automatically adjusts to the variability of
individual entries. Estimation of sparse precision matrices has also
drawn considerable attention due to its close connections to Gaussian
graphical model selection. See \citet{yl07}, \citet{y10},
\citet{rwry} and \citet{cll11}. The optimal rate
of convergence for estimating
sparse inverse covariance matrices was established in \citet{clz11}.

The rest of the paper is organized as follows. Section~\ref{secmeth}
presents a detailed construction of the data-driven block thresholding
estimator $\hat\Sigma$. The theoretical properties of the estimator
are investigated in Section~\ref{secadapt}. It is shown that the
estimator $\hat\Sigma$ achieves the optimal rate of convergence
simultaneously over each $\Ccal_\alpha(M_0, M)$ for all $\alpha,
M_0, M
> 0$. In addition, it is also shown that a slightly modified version of
$\hat\Sigma^{-1}$ is adaptively rate-optimal for estimating $\Sigma^{-1}$
over the collection $\Ccal_\alpha(M_0, M)$.
Simulation studies are carried out to illustrate the merits of the
proposed method, and the numerical results are presented in
Section~\ref{secnum}. Section~\ref{secdiss} discusses extension to subguassian
noise, adaptive estimation under the Frobenius norm and other related
issues. The proofs of the main results are given in Section~\ref{secproof}.

\section{Block thresholding}
\label{secmeth}

In this section we present in detail the construction of the adaptive
covariance matrix estimator.
The main strategy in the construction is to divide the sample
covariance matrix into blocks and then apply thresholding to each block
according to their sizes and dimensions. We shall explain these two
steps separately in Sections~\ref{blockssec} and~\ref{thresholdingsec}.

\subsection{Construction of blocks}
\label{blockssec}

As mentioned in the \hyperref[secint]{Introduction}, the application of block thresholding
to covariance matrix estimation requires more care than in the
conventional sequence estimation problems such as those from
nonparametric function estimation.\vadjust{\goodbreak} We begin with the blocking scheme
for a general $p\times p$ symmetric matrix. A key in our construction
is to make blocks larger for entries that are farther away from the
diagonal and take advantage of the approximately banding structure of
the covariance matrices in $\Ccal_\alpha$. Before we give a precise
description of the construction of the blocks, it is helpful to
graphically illustrate the construction in the following plot.

\begin{figure}

\includegraphics{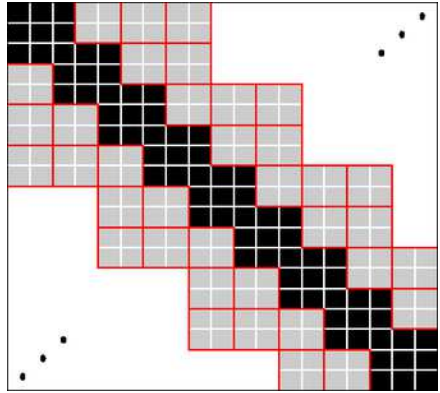}

\caption{Construction of blocks with increasing dimensions away from
the diagonal. The solid black blocks are of size $k_0\times k_0$. The
gray ones are of size $2k_0\times2k_0$.}
\label{figblock}
\end{figure}

Due to the symmetry, we shall focus only on the upper half for brevity.
We start by constructing blocks of size $k_0\times k_0$ along the
diagonal as indicated by the darkest squares in Figure~\ref{figblock}.
Note that the last block may be of a smaller size if $k_0$ is not a
divisor of $p$. Next, new blocks are created successively toward the
top right corner. We would like to increase the block sizes along the
way. To this end, we extend to the right from the diagonal blocks by
either two or one block of the same dimensions ($k_0\times k_0$) in an
alternating fashion. After this step, as exhibited in Figure~\ref
{figblock}, the odd rows of blocks will have three $k_0\times k_0$
blocks, and the even rows will have two $k_0\times k_0$ in the upper
half. Next, the size of new blocks is doubled to $2k_0\times2k_0$.
Similarly to before, the last block may be of smaller size if $2k_0$ is
not a divisor of $p$, and for the most part, we shall neglect such a
caveat hereafter for brevity. The same procedure is then followed. We
extend to the right again by three or two blocks of the size
$2k_0\times2k_0$. Afterwards, the block size is again enlarged to
$2^2k_0\times2^2k_0$ and we extend to the right by three or two blocks
of size $2^2k_0\times2^2k_0$. This procedure will continue until the
whole upper half of the $p\times p$ matrix is covered. For the lower
half, the same construction is followed to yield a symmetric blocking
of the whole matrix.

The initial block size $k_0$ can take any value as long as $k_0\asymp
\log p$. In particular, we can take $k_0=\lfloor\log p \rfloor$. The
specific choice of $k_0$ does not impact the rate of convergence, but
in practice it may be beneficial sometimes to use a value different
from $\lfloor\log p \rfloor$. In what follows, we shall keep using
$k_0$ for the sake of generality.

For notational purposes, hereafter we shall refer to the collection of
index sets for the blocks created in this fashion as $\Bcal=\{
B_1,\ldots
,B_N\}$ where $B_k=I_k\times J_k$ for some subintervals $I_k,
J_k\subset\{1,\ldots, p\}$. It is clear that $\Bcal$ forms a partition
of $\{1,2,\ldots, p\}^2$, that is,
\[
B_{k_{1}}\cap B_{k_{2}}=\varnothing\qquad \mbox{if $k_1\neq
k_2$}\quad \mbox{and} \quad B_1\cup B_2\cup\cdots\cup
B_N=\{1,2,\ldots, p\}^2.
\]
For a $p\times p$ matrix $A=(a_{ij})_{1\le i,j\le p}$ and an index set
$B=I\times J\in\Bcal$, we shall also write $A_{B}=(a_{ij})_{i\in I,
j\in J}$, a $|I|\times|J|$ submatrix of $A$. Hence $A$ is uniquely
determined by $\{A_B\dvtx B\in\Bcal\}$ and the partition $\Bcal$. With
slight abuse of notation, we shall also refer to an index set $B$ as a
block when no confusion occurs, for the sake of brevity.

Denote by $d(B)$ the dimension of $B$, that is,
\[
d(B)=\max \bigl\{\card(I),\card(J) \bigr\}.
\]
Clearly by construction, most of the blocks in $\Bcal$ are necessarily
square in that $\card(I)=\card(J)=d(B)$. The exceptions occur when the
block sizes are not divisors of $p$, which leaves the blocks along the
last row and column in rectangles rather than squares. We opt for the
more general definition of $d(B)$ to account for these rectangle blocks.

\subsection{Block thresholding}
\label{thresholdingsec}

Once the blocks are constructed, the next step is to estimate the
entries of the covariance matrix $\Sigma$, block by block, through
thresholding the corresponding blocks of the sample covariance matrix
based on the location, block size and corresponding spectral norms.

We now describe the procedure in detail. Denote by $\hat{\Sigma}$ the
block thresholding estimator, and let $B=I\times J\in\Bcal$. The
estimate of the block $\Sigma_B$ is defined as follows:
\begin{longlist}[(a)]
\item[(a)] \textit{keep the diagonal blocks}: $\hat{\Sigma}_B=\bar
{\Sigma
}_B$ if $B$ is on the diagonal, that is, $I=J$;
\item[(b)] \textit{``kill'' the large blocks}: $\hat{\Sigma}_B=\mathbf{0}$ if
$d(B)>n/\log n$;
\item[(c)] \textit{threshold the intermediate blocks}: For all other
blocks $B$, set
%
\begin{equation}
\hat{\Sigma}_B=T_{\lambda_0}(\bar{\Sigma}_B)=\bar{
\Sigma}_B\cdot \mathbb {I} \biggl( \|\bar{\Sigma}_B\| >
\lambda_{0}\sqrt{\|\bar{\Sigma }_{I\times
I}\|\|\bar{
\Sigma}_{J\times J}\|} \sqrt\frac{d(B)+\log p}{n} \biggr),\hspace*{-35pt}
\end{equation}
\end{longlist}
where $\lambda_0>0$ is a turning parameter. Our theoretical development
indicates that the resulting block thresholding estimator is optimally
rate adaptive whenever $\lambda_0$ is a sufficiently large constant. In
particular, it can be taken as fixed at $\lambda_0=6$. In practice, a
data-driven choice of $\lambda_0$ could potentially lead to further
improved finite sample performance.

It is clear from the construction that the block thresholding estimate
$\hat{\Sigma}$ is fully data-driven and does not require the knowledge
of $\alpha$. The choice of the thresholding constant $\lambda_0$ comes
from our theoretical and numerical studies. See Section~\ref{secdiss}
for more discussions on the choice of $\lambda_0$.\vadjust{\goodbreak}

We should also note that, instead of the hard thresholding operator
$T_{\lambda_0}$, more general thresholding rules can also be applied in
a similar blockwise fashion. In particular, one can use block
thresholding rules $T_{\lambda_0}(\bar{\Sigma}_B)=\bar{\Sigma}_B
\cdot
t_{\lambda_B}(\|\bar{\Sigma}_B\|)$ where
\[
\lambda_B=\lambda_{0}\sqrt{\|\bar{\Sigma}_{I\times I}
\|\|\bar{\Sigma}_{J\times J}\|} \sqrt\frac{d(B)+\log p}{n},
\]
and $t_{\lambda_B}$ is a univariate thresholding rule. Typical examples
include the soft thresholding rule $t_{\lambda_B}(z)=(|z|-\lambda_B)_+\operatorname{sgn}(z)$
and the so-called adaptive lasso rule $t_{\lambda
_B}(z)=z(1-|\lambda_B/z|^\eta)_+$ for some $\eta\ge1$, among others.
\citet{rlz09} considered entrywise universal
thresholding for estimating sparse covariance matrix. In particular,
they investigate the class of univariate thresholding rules $t_{\lambda
_B}$ such that (a) $|t_{\lambda_0}(z)|\le|z|$; (b) $t_{\lambda
_B}(z)=0$ if $|z|\le\lambda_B$; and (c) $|t_{\lambda_B}(z)-z|\le
\lambda_B$. Although we will focus on the hard thresholding rule in the
present paper for brevity, all the theoretical results developed here
apply to the more general class of block thresholding rules as well.

\section{Adaptivity}
\label{secadapt}

We now study the properties of the proposed block thresholding
estimator $\hat\Sigma$ and show that the estimator simultaneously
achieves the minimax optimal rate of convergence over the full range of
$\Ccal_\alpha$ for all $\alpha>0$.
More specifically, we have the following result.

\begin{theorem}
\label{thmain}
Let $\hat{\Sigma}$ be the block thresholding estimator of $\Sigma$ as
defined in the Section~\ref{secmeth}. Then
%
\begin{equation}
\label{optimaladaptive} \sup_{\Sigma\in\Ccal_\alpha(M_0, M)}\bE\|\hat{\Sigma}-\Sigma
\|^2\le C\min\biggl\{n^{-{2\alpha}/{(2\alpha+1)}}+\frac{\log p}{n},
\frac
{p}{n} \biggr\}
\end{equation}
for all $\alpha>0$, where $C$ is a positive constant not depending on
$n$ and $p$.
\end{theorem}
Comparing with the minimax rate of convergence given in \citet
{czz10}, this shows that the block thresholding estimator $\hat
\Sigma$ is optimally rate adaptive over $\Ccal_\alpha$ for all
$\alpha>0$.

\begin{remark}
The block thresholding estimator $\hat\Sigma$ is positive definite
with high probability, but it is not guaranteed to be positive
definite. A simple additional step, as was done in \citet{cz11},
can make the final estimator positive semi-definite and still achieve
the optimal rate of convergence. Write the eigen-decomposition of $\hat
{\Sigma}$ as $\hat{\Sigma}=\sum_{i=1}^{p}\hat{\lambda}_{i}v_{i}v_{i}^{T}$,
where $\hat{\lambda}_{i}$'s and $v_i$'s are, respectively, the
eigenvalues and eigenvectors of $\hat{\Sigma}$. Let $\hat{\lambda
}_{i}^{+}=\max(\hat\lambda_i, 0)$ be the positive part of~$\hat
{\lambda
}_{i}$, and define
\[
\hat{\Sigma}^{+}=\sum_{i=1}^{p}
\hat{\lambda}_{i}^{+}v_{i} v_{i}^{T}.
\]
Then $\hat{\Sigma}^{+}$ is positive semi-definite, and it can be shown
easily that $\hat{\Sigma}^{+}$ attains the same rate as $\hat\Sigma$.
See \citet{cz11} for further details. If a strictly positive
definite estimator is desired, one can also set $\hat{\lambda
}_{i}^{+}=\max(\hat\lambda_i, \varepsilon_n)$ for some small positive
value $\varepsilon_n$, say $\varepsilon_n=O(\log p/n)$, and the resulting
estimator $\hat{\Sigma}^{+}$ is then positive definite and attains the
optimal rate of convergence.
\end{remark}

The inverse of the covariance matrix, $\Omega:=\Sigma^{-1}$, is of
significant interest in many applications. An adaptive estimator of
$\Omega$ can also be constructed based on our proposed block
thresholding estimator. To this end, let $\hat{\Sigma}=\hat{U}\hat
{D}\hat{U}^\mathsf{T}$ be its eigen-decomposition, that is, $\hat{U}$ is
an orthogonal matrix, and $\hat{D}$ is a diagonal matrix. We propose to
estimate $\Omega$ by
\[
\hat{\Omega}=\hat{U} \operatorname{diag} \bigl(\min \bigl\{\hat{d}_{ii}^{-1},
n \bigr\} \bigr)\hat{U}^\mathsf{T},
\]
where $\hat{d}_{ii}$ is the $i$th diagonal element of $\hat{D}$. The
truncation of $\hat d_{ii}^{-1}$ is needed to deal with the case where
$\hat\Sigma$ is near singular.
The result presented above regarding $\hat{\Sigma}$ can be used to show
that $\hat{\Omega}$ adaptively achieves the optimal rate of convergence
for estimating $\Omega$.

\begin{theorem}
\label{thmain1}
Let $\hat{\Omega}$ be defined as above. Then
\[
\sup_{\Sigma\in\Ccal_\alpha}\bE\|\hat{\Omega}-\Omega\|^2\le C\min \biggl
\{n^{-{2\alpha}/{(2\alpha+1)}}+\frac{\log p}{n}, \frac{p}{n} \biggr\}
\]
for all $\alpha>0$, where $C>0$ is a constant not depending on $n$ and $p$.
\end{theorem}

The proof of the adaptivity results is somewhat involved and requires
some new technical tools. The main ideas in the theoretical analysis
can be summarized as follows:
\begin{itemize}
\item The different $\hat{\Sigma}-\Sigma$ can be decomposed into a sum
of matrices such that each matrix in the sum only consists of blocks in
$\mathcal{B}$ that are of the same size. The individual components in the
sum are then bounded separately according to their block sizes.

\item Although the sample covariance matrix $\bar{\Sigma}$ is not a
reliable estimator of $\Sigma$, its submatrix, $\bar{\Sigma}_B$, could
still be a good estimate of $\Sigma_B$. This is made precise through a
concentration inequality.

\item The analysis on the whole matrix is reduced to the analysis of a
matrix of much smaller dimensions, whose entries are the spectral norms
of the blocks, through the application of a so-called norm compression
inequality.

\item With high probability, large blocks in $\{{\bar\Sigma}_B\dvtx
B\in
\Bcal\}$, which correspond to negligible parts of the true covariance
matrix $\Sigma$, are all shrunk to zero because by construction they
are necessarily far away from the diagonal.
\end{itemize}
We shall elaborate below these main ideas in our analysis and introduce
some useful technical tools. The detailed proof is relegated to
Section~\ref{secproof}.

\subsection{Main strategy}

Recall that $\Bcal$ is the collection of blocks created using the
procedure in Section~\ref{blockssec}, and it forms a partition of $\{
1,2,\ldots,p\}^2$. We analyze the error $\hat{\Sigma}-\Sigma$ by first
decomposing it into a sum of matrices such that each matrix in the sum
only consists of blocks in $\mathcal{B}$ that are of the same size. More
precisely, for a $p\times p$ matrix $A$, define $S(A; l)$ to be a
$p\times p$ matrix whose $(i,j)$ entry equals that of $A$ if $(i,j)$
belongs to a block of dimension $2^{l-1}k_0$, and zero otherwise. In
other words,
\[
S(A,l)=\sum_{B\in\Bcal:d(B)=2^{l-1}k_0}A\circ\ind \bigl((i,j)\in B
\bigr)_{1\le i,j\le p}.
\]
With this notation, $\hat{\Sigma}-\Sigma$ is decomposed as
\[
\hat{\Sigma}-\Sigma=S(\hat{\Sigma}-\Sigma,1)+S(\hat{\Sigma }-\Sigma,2)+\cdots.
\]
This decomposition into the sum of blocks of different sizes is
illustrated in Figure~\ref{figdiv} below.

\begin{figure}
\centering
\begin{tabular}{@{}cc@{}}

\includegraphics{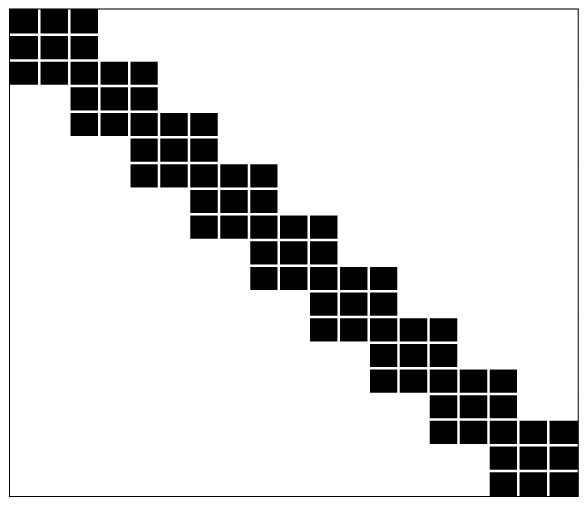}
 & \includegraphics{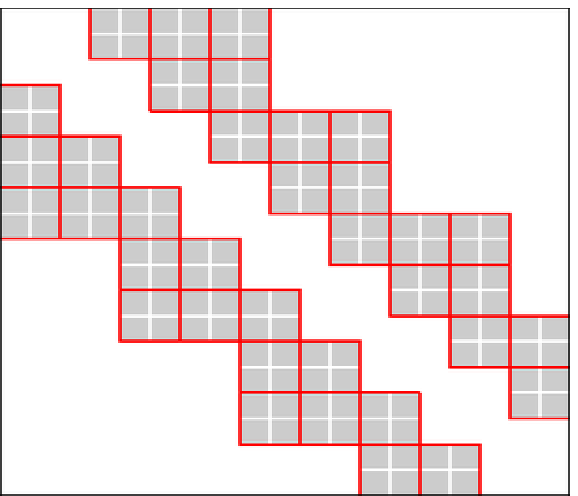}\\
\footnotesize{(a) $S(\cdot,1)$} & \footnotesize{(b) $S(\cdot,2)$}\vspace*{-3pt}
\end{tabular}
\caption{Decompose a matrix into the sum of matrices of different
block sizes: $S(\cdot,1)$ on the left and $S(\cdot,2)$ on the right.
All entries in the unshaded area are zero.}
\label{figdiv}
\end{figure}

\begin{figure}
\centering
\begin{tabular}{@{}cc@{}}

\includegraphics{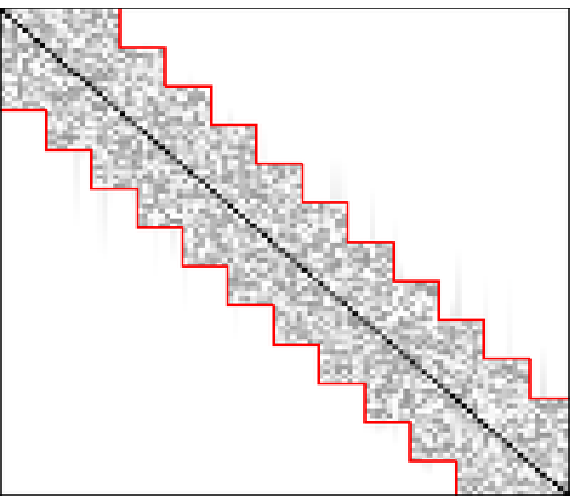}
 & \includegraphics{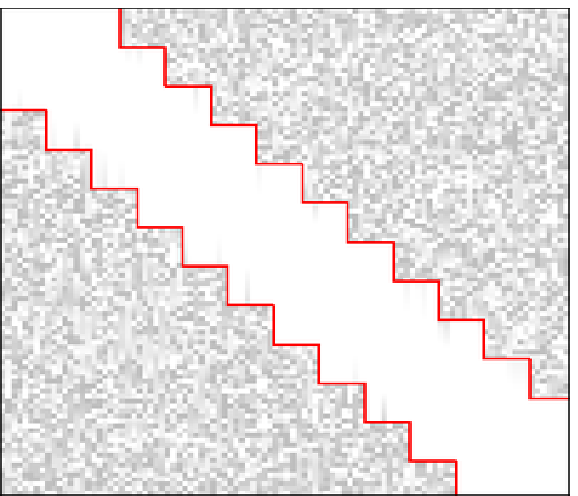}\\
\footnotesize{(a) Small blocks} & \footnotesize{(b) Large blocks}
\end{tabular}
\caption{Small blocks and large blocks are treated separately. Small
blocks are necessarily close to the diagonal and large blocks are away
from the diagonal.}
\label{figsmalllarge}
\end{figure}

We shall first separate the blocks into two groups, one for big blocks
and another for small blocks.
See Figure~\ref{figsmalllarge} for an illustration.
By the triangle inequality, for any $L\ge1$,
%
\begin{equation}
\label{eqsmalllarge} \|\hat{\Sigma}-\Sigma\|\le\sum
_{l\le L}\bigl\|S(\hat{\Sigma}-\Sigma,l)\bigr\| + \biggl\llVert\sum
_{l>L}S(\hat{\Sigma}-\Sigma,l) \biggr\rrVert.
\end{equation}
The errors on the big blocks will be bounded as a whole, and the errors
on the small blocks will be bounded separately according to block
sizes. With a careful choice of the cutoff value $L$, it can be shown
that there exists a constant $c>0$ not depending on $n$ and $p$ such
that for any $\alpha>0$ and $\Sigma\in\Ccal_\alpha$,
%
\begin{equation}
\label{eqsmall} \bE \biggl(\sum_{l\le L}\bigl\|S(\hat{
\Sigma}-\Sigma,l)\bigr\| \biggr)^2=c\min \biggl\{ n^{-{2\alpha
}/{(2\alpha+1)}}+
\frac{\log p}{n}, \frac{p}{n} \biggr\},\vadjust{\goodbreak}
\end{equation}
and
%
\begin{equation}
\label{eqlarge} \bE \biggl\llVert\sum_{l>L}S(\hat{
\Sigma}-\Sigma,l) \biggr\rrVert^2=c\min \biggl\{ n^{-{2\alpha
}/{(2\alpha+1)}}+
\frac{\log p}{n}, \frac{p}{n} \biggr\},
\end{equation}
which then implies Theorem~\ref{thmain} because
\[
\bE\|\hat{\Sigma}-\Sigma\|^2\le2\bE \biggl(\sum
_{l\le L}\bigl\|S(\hat{\Sigma}-\Sigma,l)\bigr\| \biggr)^2 + 2
\bE \biggl\llVert\sum_{l>L}S(\hat{\Sigma}-\Sigma,l)
\biggr\rrVert^2.
\]

The choice of the cutoff value $L$ depends on $p$ and $n$ and different
approaches are taken to establish~(\ref{eqsmall}) and (\ref
{eqlarge}). In both cases, a key technical tool we shall use is a
concentration inequality on the deviation of a block of the sample
covariance matrix from its counterpart of the true covariance matrix,
which we now describe.

\subsection{Concentration inequality}
\label{secconcentration}

The rationale behind our block thresholding approach is that although
the sample covariance matrix $\bar{\Sigma}$ is not a reliable estimator
of $\Sigma$, its submatrix, $\bar{\Sigma}_B$, could still be a good
estimate of~$\Sigma_B$. This observation is formalized in the
following theorem.

\begin{theorem}
\label{thunion}
There exists an absolute constant $c_0>0$ such that for all $t>1$,
\begin{eqnarray*}
&&\bP \biggl(\bigcap_{B=I\times J\in\Bcal} \biggl\{\|\bar{\Sigma
}_B-\Sigma_B\| < c_0t \sqrt{\|{
\Sigma}_{I\times I}\|\|{\Sigma}_{J\times J}\|} \sqrt\frac{d(B)+\log
p}{n} \biggr
\} \biggr)\\
&&\qquad \ge1-p^{-(6t^2-2)}.
\end{eqnarray*}
In particular, we can take $c_0=5.44$.
\end{theorem}

Theorem~\ref{thunion} enables one to bound the estimation error $\hat
{\Sigma}-\Sigma$ block by block. Note that larger blocks are
necessarily far away from the diagonal by construction. For bandable
matrices, this means that larger blocks are necessarily small in the
spectral norm. From Theorem~\ref{thunion}, if $\lambda_0>c_0$, with
overwhelming probability,
\begin{eqnarray*}
\|\bar{\Sigma}_B\|&\le& \|\Sigma_B\|+c_0
\sqrt{\|{\Sigma}_{I\times
I}\|\| {\Sigma}_{J\times J}\|} \sqrt
\frac{d(B)+\log p}{n}
\\
&<&\lambda_0\sqrt{\|{\Sigma}_{I\times I}\|\|{
\Sigma}_{J\times J}\|} \sqrt\frac{d(B)+\log p}{n}
\end{eqnarray*}
for blocks with sufficiently large sizes. As we shall show in
Section~\ref{secproof}, $\|{\Sigma}_{I\times I}\|$ and $\|{\Sigma
}_{J\times
J}\|$ in the above inequality can be replaced by their respective
sample counterparts. This observation suggests that larger blocks are
shrunken to zero with our proposed block thresholding procedure, which
is essential in establishing~(\ref{eqlarge}).

The treatment of smaller blocks is more complicated. In light of
Theorem~\ref{thunion}, blocks of smaller sizes can be estimated well,
that is, $\bar{\Sigma}_B$ is close to $\Sigma_B$ for $B$ of smaller
sizes. To translate the closeness in such a blockwise fashion into the
closeness in terms of the whole covariance matrix, we need a simple yet
useful result based on a matrix norm compression transform.

\subsection{Norm compression inequality}
\label{seccomp}

We shall now present a so-called norm compression inequality which is
particularly useful for analyzing the properties of the block
thresholding estimators.
We begin by introducing a matrix norm compression transform.

Let $A$ be a $p\times p$ symmetric matrix, and let $p_1, \ldots, p_G$ be
positive integers such that $p_1+\cdots+p_G=p$. The matrix $A$ can then
be partitioned
in a block form as
\[
A=\pmatrix{
A_{11}&A_{12}&
\ldots&A_{1G}
\vspace*{2pt}\cr
A_{21}& A_{22}& \ldots& A_{2G}
\vspace*{2pt}\cr
\vdots& \vdots& \ddots& \vdots
\vspace*{2pt}\cr
A_{G1}& A_{G2}&\ldots&A_{GG}},
\]
where $A_{ij}$ is a $p_i\times p_j$ submatrix. We shall call such a
partition of the matrix $A$ a regular partition and the blocks $A_{ij}$
regular blocks. Denote by $\Ncal\dvtx \bR^{p\times p}\mapsto\bR^{G\times
G}$ a norm compression transform
\[
A\mapsto\Ncal(A;p_1,\ldots,p_G)=\pmatrix{
\|A_{11}\|&\|A_{12}\|&\ldots&
\|A_{1G}\|
\vspace*{2pt}\cr
\|A_{21}\| & \|A_{22}\|& \ldots& \|A_{2G}\|
\vspace*{2pt}\cr
\vdots& \vdots& \ddots& \vdots
\vspace*{2pt}\cr
\| A_{G1}\|& \|A_{G2}\|&\ldots&\|A_{GG}\|}.
\]
The following theorem shows that such a norm compression transform does
not decrease the matrix norm.\vadjust{\goodbreak}

\begin{theorem}[(Norm compression inequality)]
\label{thcomp}
For any $p\times p$ matrix $A$ and block sizes $p_1,p_2,\ldots, p_G$
such that $p_1+\cdots+p_G=p$,
\[
\|A\|\le\bigl\|\Ncal(A;p_1,\ldots,p_G)\bigr\|.
\]
\end{theorem}

Together with Theorems~\ref{thunion} and~\ref{thcomp} provides a
very useful tool for bounding $S(\hat{\Sigma}-\Sigma,l)$. Note first
that Theorem~\ref{thcomp} only applies to a regular partition, that
is, the divisions of the rows and columns are the same. It is clear
that $S(\cdot,1)$ corresponds to regular blocks of size $k_0\times k_0$
with the possible exception of the last row and column which can be of
a different size, that is, $p_1=p_2=\cdots=k_0$. Hence, Theorem~\ref
{thcomp} can be directly applied. However, this is no longer the case
when $l>1$.

To take advantage of Theorem~\ref{thcomp}, a new blocking scheme is
needed for $S(\cdot,l)$. Consider the case when $l=2$. It is clear that
$S(l,2)$ does not form a regular blocking. But we can form new blocks
with $p_1=p_2=\cdots=k_0$, that is, half the size of the original
blocks in $S(\cdot, 2)$. Denote by the collection of the new blocks
$\Bcal'$. It is clear that under this new blocking scheme, each block
$B$ of size $2k_0$ consists of four elements from $\Bcal'$. Thus
\[
S(A,2)=\mathop{\sum_{B\in\Bcal}}_{d(B)=2k_0}A\circ
\ind \bigl((i,j)\in B \bigr)=\mathop{\mathop{\sum_{B'\in\Bcal'}}_{ \exists B\in
\Bcal{\mathrm{\ such\ that}\ }\,d(B)=2k_0}}_{\mathrm{and\ }B'\subset B}A
\circ\ind \bigl((i,j)\in B' \bigr).
\]

Applying Theorem~\ref{thcomp} to the regular blocks $\Bcal'$ yields
\[
\bigl\|S(A,2)\bigr\|\le \bigl\llVert\Ncal \bigl(S(A,2);k_0,
\ldots,k_0 \bigr) \bigr\rrVert,%
\]
which can be further bounded by
\[
\bigl\llVert\Ncal \bigl(S(A,2);k_0,\ldots,k_0 \bigr)
\bigr\rrVert_{\ell_1},
\]
where $\|\cdot\|_{\ell_1}$ stands for the matrix $\ell_1$ norm. Observe
that each row or column of $\Ncal(S(A,2);k_0,\ldots,k_0 )$
has at most 12 nonzero entries, and each entry is bounded by
\[
\mathop{\mathop{\max_{B'\in\Bcal'}}_{ \exists B\in\Bcal{\mathrm{\ such\ that\ }
}\,d(B)=2k_0}}_{\mathrm{and\ }B'\subset B}\|A_{B'}\|\le\mathop{
\max_{B\in
\Bcal}}_{d(B)=2k_0}\|A_B\|
\]
because $B'\subset B$ implies $\|A_{B'}\|\le\|A_B\|$. This property
suggests that $\|S(\hat{\Sigma}-\Sigma,l)\|$ can be controlled in a
block-by-block fashion. This can be done using the concentration
inequalities given in Section~\ref{secconcentration}.

The case when $l>2$ can be treated similarly. Let
$p_{2j-1}=(2^{l-1}-3)k_0$ and $p_{2j}=3k_0$ for $j=1,2,\ldots.$ It is
not hard to see that each block $B$ in $\Bcal$ of size $2^{l-1}k_0$
occupies up to four blocks in this regular blocking. And following the
same argument as before, we can derive bounds for $S(A,l)$.
%

The detailed proofs of Theorems~\ref{thmain} and~\ref{thmain1} are
given in Section~\ref{secproof}.

\section{Numerical results}\label{secnum}

The block thresholding estimator $\hat\Sigma$ proposed in
Section~\ref
{secmeth} is easy to implement. In this section we turn to the
numerical performance of the estimator. The simulation study further
illustrates the merits of the proposed block thresholding estimator.
The performance is relatively insensitive to the choice of $k_0$, and
we shall focus on $k_0=\lfloor\log p\rfloor$ throughout this section
for brevity.

We consider two different sets of covariance matrices. The setting of
our first set of numerical experiments is similar to those from
\citet{czz10}. Specifically, the true covariance matrix $\Sigma
$ is of the form
\[
\sigma_{ij}=\cases{ %
1,&\quad  $1\le
i=j\le p$,
\vspace*{2pt}\cr
\rho|i-j|^{-2}u_{ij},& \quad $1\le i\neq j\le p,$ }
\]
where the value of $\rho$ is set to be $0.6$ to ensure positive
definiteness of all covariance matrices, and $u_{ij}=u_{ji}$ are
independently sampled from a uniform distribution between $0$ and $1$.

The second settings are slightly more complicated, and the covariance
matrix $\Sigma$ is randomly generated as follows. We first simulate a
symmetric matrix $A=(a_{ij})$ whose diagonal entries are zero and
off-diagonal entries $a_{ij}$ ($i<j$) are independently generated as
$a_{ij}\sim N(0, |i-j|^{-4})$. Let $\lambda_{\min}(A)$ be its smallest
eigenvalue.
The covariance matrix $\Sigma$ is then set to be $\Sigma=\max(0,
-1.1\lambda_{\min}(A)) I+A$ to ensure its positive definiteness.

\begin{figure}

\includegraphics{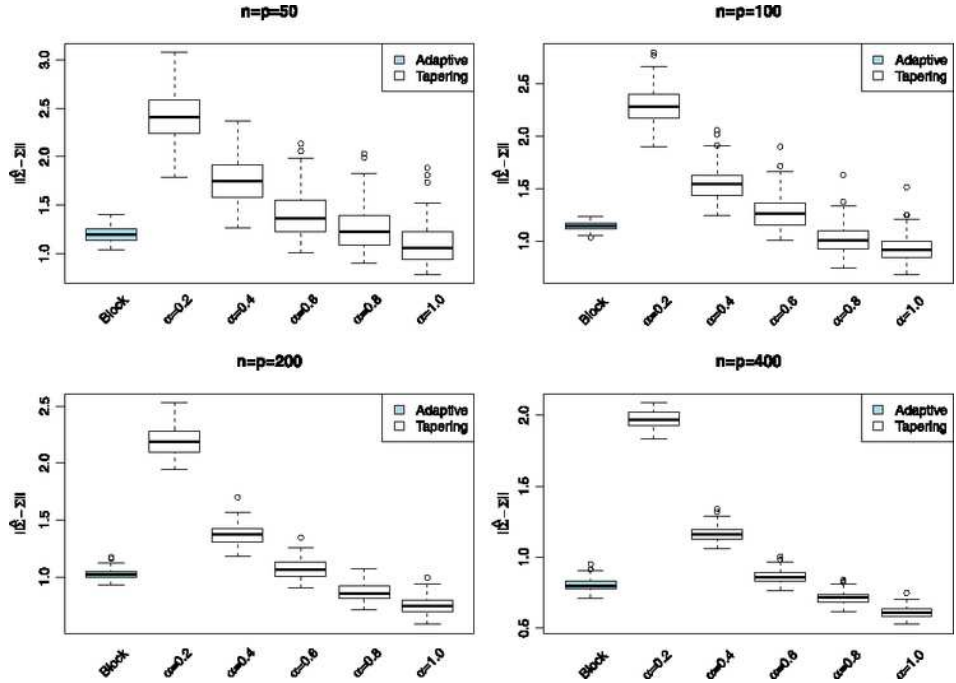}

\caption{Comparison between the tapering and adaptive block
thresholding estimators---simulation setting 1: each panel
corresponds to a particular combination of sample size $n$ and
dimension $p$. In each panel, boxplots of the estimation errors,
measured in terms of the spectral norm, are given for the block
thresholding estimator with $\lambda_0=6$ and the tapering estimator
with $\alpha=0.2$, $0.4, 0.6, 0.8$ and $1$.}
\label{figsimres}
\end{figure}

\begin{figure}

\includegraphics{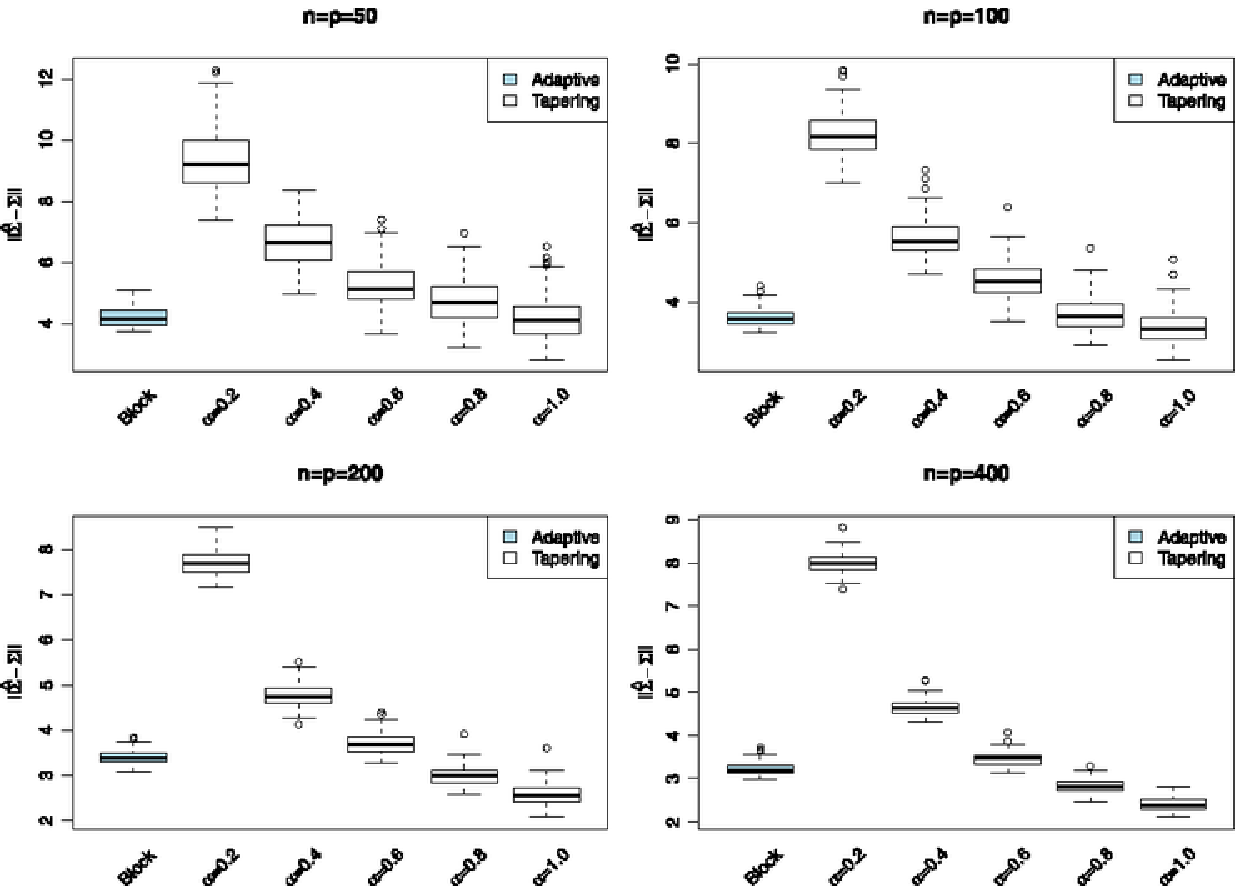}

\caption{Comparison between the tapering and adaptive block
thresholding estimators---simulation setting 2: each panel
corresponds to a particular combination of sample size $n$ and
dimension $p$. In each panel, boxplots of the estimation errors,
measured in terms of the spectral norm, are given for the block
thresholding estimator with $\lambda_0=6$ and the tapering estimator
with $\alpha=0.2$, $0.4, 0.6, 0.8$ and $1$.}
\label{fignewsimres}
\end{figure}

For each setting, four different combinations of $p$ and $n$ are
considered, $(n,p)=(50,50), (100,100), (200,200)$ and $(400,400)$, and
for each combination, 200 simulated datasets are generated. On each
simulated dataset, we apply the proposed block thresholding procedure
with $\lambda_0=6$.
For comparison purposes, we also use the banding estimator of
\citet{bl08a} and tapering estimator of \citet{czz10}
on the simulated datasets. For both estimators, a tuning parameter $k$
needs to be chosen. The two estimators perform similarly for the
similar values of $k$. For brevity, we report only the results for the
tapering estimator because it is known to be rate optimal if $k$ is
appropriately selected based on the true parameter space. It is clear
that for both our settings, $\Sigma\in\Ccal_\alpha$ with $\alpha=1$.
But such knowledge would be absent in practice. To demonstrate the
importance of knowing the true parameter space for these estimators and
consequently the necessity of an adaptive estimator such as the one
proposed here, we apply the estimators with five different values of
$\alpha$, $0.2, 0.4, 0.6, 0.8$ and $1$. We chose $k=\lfloor
n^{1/(2\alpha+1)}\rfloor$ for the tapering estimator following
\citet{czz10}.The performance of these estimators is summarized
in Figures~\ref{figsimres} and~\ref{fignewsimres} for the two
settings, respectively.

It can be seen in both settings that the numerical performance of the
tapering estimators critically depends on the specification of the
decay rate $\alpha$. Mis-specifying~$\alpha$ could lead to rather poor
performance by the tapering estimators. It is perhaps not surprising to
observe that the tapering estimator with $\alpha=1$ performed the best
among all estimators since it correctly specifies the true decay rate
and therefore, in a certain sense, made use of the information that may
not be known a priori in practice. In contrast, the proposed block
thresholding estimator yields competitive performance while not using
such information.
\section{Discussion}
\label{secdiss}

In this paper we introduced a fully data-driven covariance matrix
estimator by blockwise thresholding of the sample covariance matrix.
The estimator simultaneously attains the optimal rate of convergence
for estimating bandable covariance matrices over the full range of the
parameter spaces $\Ccal_\alpha$ for all $\alpha>0$. The estimator also
performs well numerically.

As noted in Section~\ref{thresholdingsec}, the choice of the
thresholding constant $\lambda_0=6$ is based on our theoretical and
numerical studies. Similar to wavelet thresholding in nonparametric
function estimation, in principle other choices of $\lambda_0$ can also
be used. For example, the adaptivity results on the block thresholding
estimator holds as long as
$\lambda_0\ge5.44\ (=\sqrt{24}/(1-2e^{-3}))$ where the value $5.44$
comes from the concentration inequality given in Theorem~\ref{thunion}.
Our experience suggests the performance of the block thresholding
estimator is relatively insensitive to a small change of $\lambda_0$.
However, numerically the estimator can sometimes be further improved by
using data-dependent choices of $\lambda_0$. 

Throughout the paper, we have focused on the Gaussian case for ease of
exposition and to allow for the most clear description of the block
thresholding estimator. The method and the results can also be extended
to more general subgaussian distributions. Suppose that the
distribution of the $X^{(i)}$'s is subgaussian in the sense that there
exists a constant $\sigma>0$ such that
%
\begin{equation}
\mathbb{P} \bigl\{\bigl|\mathbf{v}^\mathsf{T}(X-\mathbb{E}X)\bigr|>t \bigr\}\leq
e^{-t^{2}
/2\sigma^2}\qquad \mbox{for all }t>0\mbox{ and }\Vert\mathbf{%
v}
\Vert=1. \label{subGau}
\end{equation}
Let $\Fcal_\alpha(\sigma, M_0, M)$ denote the collection of
distributions satisfying both~(\ref{paraspace}) and~(\ref{subGau}).
Then for any given $\sigma_0>0$, the block thresholding
estimator $\hat\Sigma$ adaptively attains the optimal rate of
convergence over $\Fcal_\alpha(\sigma, M_0, M)$ for all~$\alpha$,
$M_0, M
> 0$ and $0<\sigma\le\sigma_0$ whenever $\lambda_0$ is chosen
sufficiently large.

In this paper we have focused on estimation under the spectral norm.
The block thresholding procedure, however, can be naturally extended to
achieve adaption under other matrix norms. Consider, for example, the
Frobenius norm. In this case, it is natural and also necessary to
threshold the blocks based on their respective Frobenius norms instead
of the spectral norms. Then following a similar argument as before, it
can be shown that this Frobenius norm based block thresholding
estimator can adaptively achieve the minimax rate of convergence over
every $\Ccal_\alpha$ for all $\alpha>0$. It should also be noted that
adaptive estimation under the Frobenius norm is a much easier problem
because the squared Frobenius norm is entrywise decomposable, and the
matrix can then be estimated well row by row or column by column. For
example, applying a suitable block thresholding procedure for sequence
estimation to the sample covariance matrix, row-by-row would also lead
to an adaptive covariance matrix estimator.

The block thresholding approach can also be used for estimating sparse
covariance matrices. A major difference in this case from that of
estimating bandable covariance matrices is that the block sizes cannot
be too large.
With suitable choices of the block size and thresholding level, a fully
data-driven block thresholding estimator can be shown to be
rate-optimal for estimating sparse covariance matrices. We shall report
the details of these results elsewhere.


\section{Proofs}
\label{secproof}

In this section we shall first prove Theorems~\ref{thunion} and \ref
{thcomp} and then prove the main results, Theorems~\ref{thmain} and
\ref{thmain1}. The proofs of some additional technical lemmas are
given at the end of the section.

\subsection{\texorpdfstring{Proof of Theorem \protect\ref{thunion}}{Proof of Theorem 3.3}}
The proof relies the following lemmas.

\begin{lemma}
\label{lecorr}
Let $A$ be a $2\times2$ random matrix following the Wishart
distribution $W(n,A_0)$ where
\[
A_0=\pmatrix{
1&\rho
\vspace*{2pt}\cr
\rho& 1 }.
\]
Then
\[
\bP\bigl(|A_{12}-\rho|\ge x\bigr)\le2\bP\bigl(|W_n-n|\ge{nx} \bigr),
\]
where $W_n\sim\chi^2_n$.
\end{lemma}

\begin{pf} Let $Z=(Z_1,Z_2)^\mathsf{T}\sim N(0,A_0)$ and
$Z^{(1)},\ldots, Z^{(n)}$ be $n$ independent copies of $Z$. Let
\[
S=\frac{1}{n}\sum_{i=1}^n
Z^{(i)} \bigl(Z^{(i)} \bigr)^\mathsf{T}
\]
be its sample covariance matrix. It is clear that $S=_d A$. Hence
\[
\bP\bigl(|A_{12}-\rho|\ge x\bigr)=\bP\bigl(|S_{12}-\rho|\ge x \bigr).\vadjust{\goodbreak}
\]
Note that
\[
S_{12}-\rho=\frac{1}{4} \Biggl(\frac{1}{n}\sum
_{i=1}^n \bigl( \bigl(Z_1^{(i)}+Z_2^{(i)}
\bigr)^2-2(1+\rho) \bigr)-\frac{1}{n}\sum
_{i=1}^n \bigl( \bigl(Z_1^{(i)}-Z_2^{(i)}
\bigr)^2-2(1-\rho) \bigr) \Biggr).
\]
Therefore,
\begin{eqnarray*}
&&\bP\bigl(|S_{12}-\rho|\ge x \bigr)\\
&&\qquad\le\bP \Biggl( \Biggl\vert \frac{1}{n}\sum
_{i=1}^n \bigl( \bigl(Z_1^{(i)}+Z_2^{(i)}
\bigr)^2-2(1+\rho) \bigr) \Biggr\vert\ge2(1+\rho)x \Biggr)
\\
&&\qquad\quad{}+\bP \Biggl( \Biggl\vert\frac{1}{n}\sum_{i=1}^n
\bigl( \bigl(Z_1^{(i)}-Z_2^{(i)}
\bigr)^2-2(1-\rho) \bigr) \Biggr\vert\ge2(1-\rho)x \Biggr).
\end{eqnarray*}
Observe that
\begin{eqnarray*}
&&\bP \Biggl(\Biggl \vert\frac{1}{n}\sum_{i=1}^n
\bigl( \bigl(Z_1^{(i)}+Z_2^{(i)}
\bigr)^2-2(1+\rho) \bigr) \Biggr\vert\ge2(1+\rho)x \Biggr)
\\
&&\qquad=\bP \Biggl( \Biggl\vert\sum_{i=1}^n
\frac{ (Z_1^{(i)}+Z_2^{(i)}
)^2}{ 2(1+\rho)}-n \Biggr\vert\ge x \Biggr)
\\
&&\qquad=\bP\bigl(|W_n-n|\ge x\bigr).
\end{eqnarray*}
Similarly,
\[
\bP \Biggl( \Biggl\vert\frac{1}{n}\sum_{i=1}^n
\bigl( \bigl(Z_1^{(i)}-Z_2^{(i)}
\bigr)^2-2(1-\rho) \bigr) \Biggr\vert\ge2(1-\rho)x \Biggr)=
\bP\bigl(|W_n-n|\ge x\bigr).
\]
The proof is now complete.
\end{pf}

\begin{lemma}
\label{letail}
Let $B=I\times J\subset[1,p]^2$. There exists an absolute constant
$c_0>0$ such that for any $t>1$,
\[
\bP \biggl\{\|\bar{\Sigma}_B-\Sigma_B\| <
c_0t \sqrt{\|{\Sigma}_{I\times
I}\|\|{\Sigma}_{J\times J}\|}
\sqrt\frac{d(B)+\log p}{n} \biggr\} \ge1-p^{-6t^2}.
\]
In particular, we can take $c_0=5.44$.
\end{lemma}

\begin{pf} Without loss of generality, assume that
$\card
(I)=\card(J)=d(B)=d$. Let $A$ be a $d\times d$ matrix, $\bu_1,\bu_2$
and $\bv_1,\bv_2\in\Scal^{d-1}$ where $\Scal^{d-1}$ is the unit sphere
in the $d$ dimensional Euclidean space. Observe that
\begin{eqnarray*}
\bigl|\bu_1^\mathsf{T}A\bv_1\bigr|-\bigl|\bu_2^\mathsf{T}A
\bv_2\bigr|&\le& \bigl|\bu_1^\mathsf{T}A\bv_1-
\bu_2^\mathsf{T}A\bv_2\bigr|
\\
&=& \bigl\vert\bu_1^\mathsf{T} A(\bv_1-
\bv_2)+(\bu_1-\bu_2)^\mathsf{T}A
\bv_2 \bigr\vert
\\
&\le& \bigl\vert\bu_1^\mathsf{T} A(\bv_1-
\bv_2) \bigr\vert+ \bigl\vert(\bu_1- \bu_2)^\mathsf{T}A
\bv_2 \bigr\vert
\\
&\le&\|\bu_1\|\|A\|\|\bv_1-\bv_2\|+\|
\bu_1-\bu_2\|\|A\|\|\bv_2\|
\\
&=&\|A\| \bigl(\|\bv_1-\bv_2\|+\|\bu_1-
\bu_2\| \bigr),
\end{eqnarray*}
where as before, we use $\|\cdot\|$ to represent the spectral norm for
a matrix and $\ell_2$ norm for a vector. As shown by B\"or\"oczky and
Wintsche [(\citeyear{bw05}), e.g., Corollary~1.2], there exists an
$\delta$-cover
set $Q_d\subset\Scal^{d-1}$ of $\Scal^{d-1}$ such that
\[
\card(Q_d)\le\frac{c\cos\delta}{\sin^d\delta}d^{3/2}\log \bigl(1+d
\cos^2\delta \bigr)\approx c\delta^{-d}d^{3/2}
\log(1+d)
\]
for some absolute constant $c>0$. Note that
%
\begin{equation}
\|A\|=\sup_{\bu,\bv\in\Scal^{d-1}} \bu^\mathsf{T}A\bv\le\sup_{\bu
,\bv\in
Q_{d}}
\bu^\mathsf{T}A\bv+2\delta\|A\|.
\end{equation}
In other words,
%
\begin{equation}
\|A\|\le(1-2\delta)^{-1}\sup_{\bu,\bv\in Q_{d}} \bu^\mathsf{T}A\bv.
\end{equation}

Now consider $A=\bar{\Sigma}_B-\Sigma_B$. Let $X_I=(X_i\dvtx i\in
I)^\mathsf{T}$ and $X_J=(X_i\dvtx i\in J)^\mathsf{T}$. Then
\[
\bar{\Sigma}_B=\frac{1}{n}\sum_{i=1}^n
\bigl(X_I^{(i)}-\bar{X}_I \bigr)
\bigl(X_J^{(i)}-\bar{X}_J \bigr)^\mathsf{T},
\]
where
\[
\bar{X}_I=(\bar{X}_i\dvtx i\in I)^\mathsf{T}\quad \mbox{and}\quad \bar{X}_J=(\bar{X}_i\dvtx i\in J)^\mathsf{T}.
\]
Similarly, $\Sigma_B=\bE(X_I-\bE X_I)(X_J-\bE X_J)^\mathsf{T}$. Therefore,
\[
A=\frac{1}{n}\sum_{i=1}^n
\bigl(X_I^{(i)} \bigl(X_J^{(i)}
\bigr)^\mathsf{T}-\bE X_IX_J^\mathsf{T} \bigr)-
\bigl(\bar{X}_I\bar{X}_J^\mathsf{T}-\bE
X_I\bE X_J^\mathsf{T} \bigr).
\]
Clearly the distributional properties of $A$ are invariant to the mean
of $X$. We shall therefore assume without loss of generality that $\bE
X=0$ in the rest of the proof.

For any fixed $\bu, \bv\in\Scal^{d-1}$, we have
\[
\bu^\mathsf{T}A\bv=\frac{1}{n}\sum_{i=1}^n
\bigl(Y_1^{(i)}Y_2^{(i)}-\bE
Y_1Y_2 \bigr)-\bar{Y}_1
\bar{Y}_2,
\]
where $Y_1=\bu^\mathsf{T}X_I$, $Y_2=\bv^\mathsf{T}X_J$, and similarly,
$Y_1^{(i)}=\bu^\mathsf{T}X_I^{(i)}$, $Y_2^{(i)}=\bv^\mathsf{T}X_J^{(i)}$. It
is not hard to see that
\[
\pmatrix{Y_1
\vspace*{2pt}\cr
Y_2 }
\sim N\lleft(\mathbf{0}, \pmatrix{ \bu^\mathsf{T}\Sigma_{I\times I}\bu&
\bu^\mathsf{T}\Sigma_{I\times J}\bv
\vspace*{2pt}\cr
\bv^\mathsf{T}\Sigma_{J\times I}\bu& \bv^\mathsf{T}
\Sigma_{J\times J}\bv } \rright),
\]
and $\bu^\mathsf{T}A\bv$ is simply the difference between the sample and
population covariance of $Y_1$ and $Y_2$. We now appeal to the
following lemma:\vadjust{\goodbreak}

Applying Lemma~\ref{lecorr}, we obtain
\[
\bP \bigl\{ \bigl\vert\bu^\mathsf{T}A\bv \bigr\vert\ge x \bigr\}\le 2\bP
\biggl(|W_n-n|\ge\frac{nx}{ (
(\bu^\mathsf{T}\Sigma_{I\times I}\bu)(\bv^\mathsf{T}\Sigma_{J\times
J}\bv
) )^{1/2}} \biggr),
\]
where $W_n\sim\chi_n^2$. By the tail bound for $\chi^2$ random
variables, we have
\[
\bP \biggl(|W_n-n|\ge\frac{nx}{ (
(\bu^\mathsf{T}\Sigma_{I\times I}\bu)(\bv^\mathsf{T}\Sigma_{J\times
J}\bv
) )^{1/2}} \biggr)\le\exp \biggl(-
\frac{nx^2}{4\|\Sigma_{I\times
I}\|
\|\Sigma_{J\times J}\|} \biggr).
\]
See, for example, Lemma~1 of Laurent and Massart (\citeyear{lm00}). In summary,
\[
\bP \bigl\{ \bigl\vert\bu^\mathsf{T}A\bv \bigr\vert\ge x \bigr\}\le 2\exp \biggl(-
\frac{nx^2}{4\|\Sigma_{I\times I}\|\|\Sigma_{J\times J}\|} \biggr).
\]

Now an application of union bound then yields
\begin{eqnarray*}
\bP\bigl(\|\bar{\Sigma}_B-\Sigma_B\|\ge x\bigr)&\le& \bP \Bigl\{
\sup_{\bu
,\bv\in
Q_{d}} \bu^\mathsf{T}A\bv\ge(1-2\delta)x \Bigr\}
\\
&\le& 2\card(Q_{d})^2\exp \biggl(-\frac{n(1-2\delta)^2x^2}{4\|
\Sigma_{I\times I}\|\|\Sigma_{J\times J}\|}
\biggr)
\\
&\le& c\delta^{-2d}d^3\log^2(1+d)\exp \biggl(-
\frac{n(1-2\delta
)^2x^2}{4\|
\Sigma_{I\times I}\|\|\Sigma_{J\times J}\|} \biggr)
\end{eqnarray*}
for some constant $c>0$. In particular, taking
\[
x=c_0t \sqrt{\|{\Sigma}_{I\times I}\|\|{\Sigma}_{J\times J}
\|} \sqrt{ \frac{d+\log p}{ n}}
\]
yields
\[
\bP\bigl(\|\bar{\Sigma}_B-\Sigma_B\|\ge x\bigr)\le c
\delta^{-2d}d^3\log^2(1+d)\exp \biggl(-
\frac{c_0^2t^2}{4}(1-2\delta)^2(d+\log p) \biggr).
\]
Let $\delta=e^{-3}$ and
\[
c_0>\frac{\sqrt{24}}{1-2\delta}=5.44.
\]
Then
\[
\bP\bigl(\|\bar{\Sigma}_B-\Sigma_B\|\ge x\bigr)\le
p^{-6t^{2}}.
\]
\upqed\end{pf}

We are now in position to prove Theorem~\ref{thunion}. It is clear
that the total number of blocks can be upper bounded by $\card(\Bcal
)\le(p/k_0)^2< p^2$. It follows from the union bound and Lemma~\ref
{letail} that
\begin{eqnarray*}
&&\bP \biggl\{\bigcup_{B\in\Bcal}\|\bar{
\Sigma}_{B}- \Sigma_{B}\|\ge c_0t\bigl (\|
\Sigma_{I\times I}\|\| \Sigma_{J\times J}\| \bigr)^{1/2}
\bigl(n^{-1} \bigl(d({B})+\log p \bigr) \bigr)^{1/2} \biggr\}
\\
&&\qquad\le\sum_{B\in\Bcal}\bP \bigl\{\|\bar{
\Sigma}_{B}-\Sigma_{B}\| \ge c_0t \bigl(\|
\Sigma_{I\times I}\|\|\Sigma_{J\times J}\| \bigr)^{1/2}
\bigl(n^{-1} \bigl(d({B})+\log p \bigr) \bigr)^{1/2} \bigr\}
\\
&&\qquad\le p^{-6t^2+2}.
\end{eqnarray*}

\subsection{\texorpdfstring{Proof of Theorem \protect\ref{thcomp}}{Proof of Theorem 3.4}}
Denote by $\bu,\bv$ the left and right singular vectors corresponding
to the leading singular value of $A$, that is, $\bu^\mathsf{T}A\bv=\|A\|$.
Let $\bu=(\bu_1,\ldots,\bu_G)^\mathsf{T}$ and $\bv=(\bv_1,\ldots
,\bv_G)^\mathsf{T}$ be partitioned in the same fashion as $X$, for example,
$\bu_g, \bv_g\in\bR^{p_g}$. Denote by $\bu_\ast=(\|\bu_1\|
,\ldots,\|\bu_G\|)^\mathsf{T}$ and $\bv_\ast=(\|\bv_1\|,\ldots,\|\bv_G\|)^\mathsf{T}$. It
is clear that $\|\bu_\ast\|=\|\bv_\ast\|=1$. Therefore,
\begin{eqnarray*}
\bigl\|\Ncal(A)\bigr\| &\ge& \bu_\ast^\mathsf{T}\Ncal(A)
\bv_\ast=\sum_{j,k=1}^G \|
\bu_j\|\|\bv_k\|\|\Sigma_{jk}\|
\\
&\ge& \sum_{j,k=1}^G \bu_j^\mathsf{T}
\Sigma_{jk}\bv_k=\bu^\mathsf{T}\Sigma\bv =\|A\|.
\end{eqnarray*}

\subsection{\texorpdfstring{Proof of Theorem \protect\ref{thmain}}{Proof of Theorem 3.1}}

With the technical tools provided by Theorems~\ref{thunion} and \ref
{thcomp}, we now show that $\hat{\Sigma}$ is an adaptive estimator of
$\Sigma$ as claimed by Theorem~\ref{thmain}. We begin by establishing
formal error bounds on the blocks using the technical tools introduced earlier.

\subsubsection{Large blocks}
\label{seclarge}

First treat the larger blocks. When $\Sigma\in\Ccal_\alpha$, large
blocks can all be shrunk to zero because they necessarily occur far
away from the diagonal and therefore are small in spectral norm. More
precisely, we have:

\begin{lemma}
\label{leblocknorm}
For any $B\in\Bcal$ with $d(B)\ge2k_0$,
\[
\|\Sigma_B\|\le Md(B)^{-\alpha}.
\]
\end{lemma}

Together with Theorem~\ref{thunion}, this suggests that
\begin{eqnarray*}
\|\bar{\Sigma}_B\|&\le& \|\bar{\Sigma}_B-
\Sigma_B\|+\|\Sigma_B\|
\\
&\le&c_0 \bigl(\|\Sigma_{I\times I}\|\|\Sigma_{J\times J}\|
\bigr)^{1/2} \bigl(n^{-1} \bigl(d({B})+\log p \bigr)
\bigr)^{1/2}+Md(B)^{-\alpha},
\end{eqnarray*}
with probability at least $1-p^{-4}$. Therefore, when
%
\begin{equation}
\label{eqdefbig} d(B)\ge c\min \biggl\{n^{{1}/{(2\alpha+1)}}, \biggl(
\frac{n}{\log
p} \biggr)^{{1} /{(2\alpha)}} \biggr\}
\end{equation}
for a large enough constant $c>0$,
%
\begin{equation}
\label{eqshrunken} \|\bar{\Sigma}_B\|<\tfrac{1}{2}(c_0+
\lambda_0) \bigl(\|\Sigma_{I\times I}\| \|\Sigma_{J\times J}\|
\bigr)^{1/2} \bigl(n^{-1} \bigl(d({B})+\log p \bigr)
\bigr)^{1/2}.
\end{equation}
The following lemma indicates that we can further replace $\|\Sigma_{I\times I}\|$ and $\|\Sigma_{J\times J}\|$ by their respective sample
counterparts.

\begin{lemma}
\label{letune}
Denote by $\Ical=\{I\dvtx I\times J\in\Bcal\}$. Then for all $I\in
\Ical$,
\[
1-\frac{\sqrt{\card(I)}+t}{\sqrt{n}}\le\frac{\|\bar{\Sigma
}_{I\times
I}\|
}{\|\Sigma_{I\times I}\|}\le1+\frac{\sqrt{\card(I)}+t}{\sqrt{n}},
\]
with probability at least $1-4p^2\exp(-t^2/2)$.
\end{lemma}

In the light of Lemma~\ref{letune},~(\ref{eqshrunken}) implies that,
with probability at least $1-2p^{-4}$, for any $B\in\Bcal$ such that
$d(B)\le n/\log n$ and~(\ref{eqdefbig}) holds,
\[
\|\bar{\Sigma}_B\|<\lambda_0 \bigl(\|\bar{
\Sigma}_{I\times I}\|\| \bar{\Sigma}_{J\times J}\| \bigr)^{1/2}
\bigl(n^{-1} \bigl(d({B})+\log p \bigr) \bigr)^{1/2},
\]
whenever $n/\log p$ is sufficiently large. In other words, with
probability at least $1-2p^{-4}$, for any $B\in\Bcal$ such that (\ref
{eqdefbig}) holds, $\hat{\Sigma}_B=\mathbf{0}$.

\subsubsection{Small blocks}
\label{secsmall}

Now consider the smaller blocks. From the discussions in Section~\ref
{seccomp}, we have
%
\begin{equation}
\label{eqncompsmall} \bigl\|S(\hat{\Sigma}-\Sigma,l)\bigr\|\le12\max_{B\in
\Bcal
\dvtx d({B})=2^{l-1}k_0}\|
\hat{\Sigma}_{B}-\Sigma_{B}\|.
\end{equation}

Observe that by the definition of $\hat{\Sigma}$,
\begin{eqnarray*}
\|\hat{\Sigma}_{B}-\Sigma_{B}\|&\le& \|\hat{
\Sigma}_{B}-\bar{\Sigma}_{B}\|+\|\bar{ \Sigma}_{B}-
\Sigma_{B}\|
\\
&\le&\lambda_0\bigl(\|\bar{\Sigma}_{I\times I}\|\|\bar{\Sigma
}_{J\times J}\| \bigr)^{1/2} \bigl(n^{-1} \bigl(d(B)+\log p
\bigr) \bigr)^{1/2}+\|\bar{\Sigma}_{B}-\Sigma_{B}
\|.
\end{eqnarray*}
By Lemma~\ref{letune}, the spectral norm of $\bar{\Sigma}_{I\times I}$
and $\bar{\Sigma}_{J\times J}$ appeared in the first term on the
rightmost-hand side can be replaced by their corresponding population
counterparts, leading to
\begin{eqnarray*}
\|\hat{\Sigma}_{B}-\Sigma_{B}\|&\le& \lambda_0\bigl(
\|{\Sigma}_{I\times I}\| \|{\Sigma}_{J\times J}\|\bigr)^{1/2}
\bigl(n^{-1} \bigl(d(B)+\log p \bigr) \bigr)^{1/2}+\|\bar{
\Sigma }_{B}-\Sigma_{B}\|
\\
&\le&\lambda_0M_0 \bigl(n^{-1} \bigl(d(B)+
\log p \bigr) \bigr)^{1/2}+\|\bar{\Sigma}_{B}-
\Sigma_{B} \|,
\end{eqnarray*}
where we used the fact that $\|{\Sigma}_{I\times I}\|, \|{\Sigma
}_{J\times J}\|\le M_0$. This can then be readily bounded, thanks to
Theorem~\ref{thunion}:%
\[
\|\hat{\Sigma}_{B}-\Sigma_{B}\|\le(\lambda_0M_0+c_0)
\bigl(n^{-1} \bigl(d(B)+\log p \bigr) \bigr)^{1/2}.
\]
Together with~(\ref{eqncompsmall}), we get
%
\begin{equation}
\label{eqsmallbd} \bigl\|S(\hat{\Sigma}-\Sigma,l)\bigr\|\le C \bigl(n^{-1}
\bigl(k_02^{l-1}+\log p \bigr) \bigr)^{1/2}.
\end{equation}

%
%

%

\subsubsection{Bounding the estimation error}
To put the bounds on both small and big blocks together, we need only
to choose an appropriate cutoff $L$ in~(\ref{eqsmalllarge}). In
particular, we take
%
\begin{equation}
\label{eqL}  L=\cases{
\bigl\lceil
\log_2 (p/k_0) \bigr\rceil,& \quad $\mbox{if }p\le
n^{1/(2\alpha+1)},$
\vspace*{2pt}\cr
\bigl\lceil\log_2 \bigl(n^{{1}/ {2\alpha
+1}}/k_0 \bigr) \bigr\rceil,&\quad $\mbox{if }\log p< n^{1/(2\alpha+1)} \mbox{ and } n^{1/(2\alpha+1)}\le p,$
\vspace*{2pt}\cr
\bigl\lceil\log_2(\log p/k_0) \bigr\rceil,&\quad $\mbox{if }n^{1/(2\alpha+1)}\le\log
p,$}\hspace*{-35pt}
\end{equation}
where $\lceil x\rceil$ stands for the smallest integer that is no less
than $x$.

\subsubsection*{Small $p$}

If $p\le n^{1/(2\alpha+1)}$, all blocks are small. From the bound
derived for small blocks, for example, equation~(\ref{eqsmallbd}), we have
\[
\|\hat{\Sigma}-\Sigma\| \le\sum_{l}\bigl\|S(\hat{
\Sigma}-\Sigma,l)\bigr\| \le C\sum_{l}
\bigl(n^{-1} \bigl(2^{l-1}k_0+\log p \bigr)
\bigr)^{1/2}
\le C(p/n)^{1/2},
\]
with probability at least $1-2p^{-4}$. Hereafter we use $C>0$ as a
generic constant that does not depend on $p$, $n$ or $\alpha$, and its
value may change at each appearance. Thus
\begin{eqnarray*}
\bE\|\hat{\Sigma}-\Sigma\|^2&=&\bE\|\hat{\Sigma}-\Sigma
\|^2{\mathbb I} \bigl\{\|\hat{\Sigma}-\Sigma\|\le C(p/n)^{1/2}
\bigr\}\\
&&{}+\bE\| \hat{\Sigma}-\Sigma\|^2{\mathbb I} \bigl\{\|\hat{
\Sigma}-\Sigma\|> C(p/n)^{1/2} \bigr\}.
\end{eqnarray*}
It now suffices to show that the second term on the right-hand side is
$O(p/n)$. By the Cauchy--Schwarz inequality,
\begin{eqnarray*}
&&\bE\|\hat{\Sigma}-\Sigma\|^2{\mathbb I} \bigl(\|\hat{\Sigma }-\Sigma
\|> C(p/n)^{1/2} \bigr)\\
&&\qquad\le \bigl(\bE\|\hat{\Sigma}-\Sigma
\|^4{\mathbb P} \bigl\{ \|\hat{\Sigma}-\Sigma\|> C(p/n)^{1/2}
\bigr\} \bigr)^{1/2}
\\
&&\qquad\le \bigl(2p^{-4}\bE\|\hat{\Sigma}-\Sigma\|^4
\bigr)^{1/2}.
\end{eqnarray*}
Observe that
\[
\bE\|\hat{\Sigma}-\Sigma\|^4\le\bE\|\hat{\Sigma}-\Sigma
\|_\mathrm{F}^4\le Cp^4/n^2,
\]
where $\|\cdot\|_\mathrm{F}$ stands for the Frobenius norm of a matrix. Thus,
\[
\bE\|\hat{\Sigma}-\Sigma\|^2{\mathbb I} \bigl\{\|\hat{\Sigma}-
\Sigma\|> C(p/n)^{1/2} \bigr\}\le Cp/n.
\]

\subsubsection*{Medium $p$}
When $\log p< n^{1/(2\alpha+1)}$ and $n^{1/(2\alpha+1)}\le p$, by the
analysis from Section~\ref{seclarge}, all large blocks will be shrunk
to zero with overwhelming probability, that is,
\[
\bP \biggl\{\sum_{l>L} S(\hat{\Sigma},l)=0 \biggr\}
\ge1-2p^{-4}.
\]
When this happens,
\[
\biggl\llVert\sum_{l>L} S(\hat{\Sigma}-\Sigma,l)
\biggr\rrVert= \biggl\llVert\sum_{l>L} S(\Sigma,l)
\biggr\rrVert\le \biggl\llVert\sum_{l>L} S(\Sigma,l)
\biggr\rrVert_{\ell_1}.
\]
Recall that $\|\cdot\|_{\ell_1}$ stands for the matrix $\ell_1$ norm,
that is, the maximum row sum of the absolute values of the entries of a
matrix. Hence,
\[
\biggl\llVert\sum_{l>L} S(\hat{\Sigma}-\Sigma,l)
\biggr\rrVert\le ML^{-\alpha
}\le Cn^{-{\alpha}/{(2\alpha+1)}}.
\]
As a result,
\begin{eqnarray*}
\bE \biggl\llVert\sum_{l>L} S(\hat{\Sigma}-\Sigma,l)
\biggr\rrVert^2&=& \bE \biggl\llVert\sum
_{l>L} S(\hat{\Sigma}-\Sigma,l) \biggr\rrVert^2{
\mathbb I} \biggl\{ \sum_{l>L}S(\hat{\Sigma},l)=0
\biggr\}
\\
&&{} + \bE \biggl\llVert\sum_{l>L} S(\hat{\Sigma}-
\Sigma,l) \biggr\rrVert^2{\mathbb I} \biggl\{\sum
_{l>L}S(\hat{\Sigma},l)\neq0 \biggr\}.
\end{eqnarray*}
It remains to show that
\[
\bE \biggl\llVert\sum_{l>L} S(\hat{\Sigma}-\Sigma,l)
\biggr\rrVert^2{\mathbb I} \biggl\{\sum_{l>L}S(
\hat{\Sigma},l)\neq0 \biggr\}=O \bigl(n^{-{2\alpha }/{(2\alpha+1)}} \bigr).
\]
By the Cauchy--Schwarz inequality,
\begin{eqnarray*}
&&\bE \biggl\{ \biggl\llVert\sum_{l>L} S(\hat{
\Sigma}-\Sigma,l) \biggr\rrVert^2{\mathbb I} \biggl(\sum
_{l>L}S(\hat{\Sigma},l)\neq0 \biggr) \biggr\}
\\
&&\qquad\le \biggl(\bE \biggl\llVert\sum_{l>L} S(\hat{
\Sigma}-\Sigma,l) \biggr\rrVert^4\bP \biggl\{\sum
_{l>L}S(\hat{\Sigma},l)\neq0 \biggr\} \biggr)^{1/2}.
\end{eqnarray*}
Observe that
\begin{eqnarray*}
\biggl\llVert\sum_{l>L} S(\hat{\Sigma}-\Sigma,l)
\biggr\rrVert^4&\le& \biggl\llVert\sum_{l>L}
S(\hat{\Sigma}-\Sigma,l) \biggr\rrVert_\mathrm{F}^4 = \biggl(
\biggl\llVert\sum_{l>L} S(\hat{\Sigma}-\Sigma,l)
\biggr\rrVert_\mathrm{F}^2 \biggr)^2
\\
&\le& \biggl( \biggl\llVert\sum_{l>L} S(\bar{
\Sigma}-\Sigma,l) \biggr\rrVert^2_\mathrm{F}+ \biggl\llVert\sum
_{l>L} S(\Sigma,l) \biggr\rrVert^2_\mathrm{F}
\biggr)^2
\\
&\le&2 \biggl( \biggl\llVert\sum_{l>L} S(\bar{
\Sigma}-\Sigma,l) \biggr\rrVert^4_\mathrm{F}+ \biggl\llVert\sum
_{l>L} S(\Sigma,l) \biggr\rrVert^4_\mathrm{F} \biggr),
\end{eqnarray*}
where the second inequality follows from the fact that $\hat{\Sigma
}=\bar
{\Sigma}$ or $\mathbf{0}$. It is not hard to see that
\[
\bE \biggl\llVert\sum_{l>L}
S(\bar{\Sigma}-\Sigma,l) \biggr\rrVert_\mathrm{F}^4 \le\bE
\llVert\bar{\Sigma}-\Sigma\rrVert_\mathrm{F}^4 \le
Cp^4/n^2.
\]
On the other hand,
\begin{eqnarray*}
\biggl\llVert\sum_{l>L} S(\Sigma,l) \biggr
\rrVert_\mathrm{F}^4&\le& \biggl(\sum_{i,j:
|i-j|>k_02^{L-1}}
\sigma_{ij}^2 \biggr)^2\le \biggl(\sum
_{i,j:
|i-j|>k_02^{L-1}} |\sigma_{ij}| \biggr)^4\\
&\le&
Cn^{-{4\alpha}/{
(2\alpha+1)}}.
\end{eqnarray*}
Therefore,
\[
\bE \biggl\llVert\sum_{l>L} S(\hat{\Sigma}-\Sigma,l)
\biggr\rrVert^4\le C \bigl({p^4}/n^2+n^{-{4\alpha}/{(2\alpha+1)}}
\bigr).
\]
Together with Theorem~\ref{thunion}, we conclude that
\[
\bE \biggl\llVert\sum_{l>L} S(\hat{\Sigma}-\Sigma,l)
\biggr\rrVert^2{\mathbb I} \biggl\{\sum_{l>L}S(
\hat{\Sigma},l)\neq0 \biggr\}\le Cn^{-1}.
\]

\subsubsection*{Large $p$}

Finally, when $p$ is very large in that $\log p>n^{1/(2\alpha+1)}$, we
can proceed in the same fashion. Following the same argument as before,
it can be shown that
\[
\bE \biggl\llVert\sum_{l>L} S(\hat{\Sigma},l) \biggr
\rrVert^2\le C \bigl(n^{-1}\log p \bigr).
\]
The smaller blocks can also be treated in a similar fashion as before.
From equation~(\ref{eqsmallbd}),
\[
\sum_{l\le L}\bigl\|S(\hat{\Sigma}-\Sigma,l)\bigr\|\le C
\bigl(n^{-1}\log p \bigr),
\]
with probability at least $1-2p^{-4}$. Thus, it can be calculated that
\[
\bE \biggl(\sum_{l\le L}\bigl\|S(\hat{\Sigma}-\Sigma,l)\bigr\|
\biggr)^2\le C \bigl(n^{-1}\log p \bigr).
\]
Combining these bounds, we conclude that
$
\bE\|\hat{\Sigma}-\Sigma\|^2\le C(n^{-1}\log p).
$
In summary,
\[
\sup_{\Sigma\in\Ccal_\alpha}\bE\|\hat{\Sigma}-\Sigma\|^2\le C\min \biggl
\{n^{-2\alpha/(2\alpha+1)}+\frac{\log p}{n},\frac{p}{n} \biggr\},
\]
for all $\alpha>0$. In other words, the block thresholding estimator
$\hat{\Sigma}$ achieves the optimal rate of convergence simultaneously
over every $\Ccal_\alpha$ for all $\alpha>0$.

\subsection{\texorpdfstring{Proof of Theorem \protect\ref{thmain1}}{Proof of Theorem 3.2}}

Observe that
\begin{eqnarray*}
\bE\|\hat{\Omega}-\Omega\|^2&=&\bE \bigl(\|\hat{\Omega}-\Omega
\|^2\ind \bigl\{ \lambda_{\min}(\hat{\Sigma})\ge
\tfrac{1}{2}\lambda_{\min}(\Sigma) \bigr\} \bigr)
\\
&&{}+\bE \bigl(\|\hat{\Omega}-\Omega\|^2\ind \bigl\{
\lambda_{\min}(\hat{\Sigma})< \tfrac{1}{2}\lambda_{\min}(
\Sigma) \bigr\} \bigr),
\end{eqnarray*}
where $\lambda_{\min}(\cdot)$ denotes the smallest eigenvalue of a
symmetric matrix. Under the event that
\[
\lambda_{\min}(\hat{\Sigma})\ge\tfrac{1}{2}\lambda_{\min}(
\Sigma),
\]
$\hat{\Sigma}$ is positive definite and $\hat{\Omega}=\hat{\Sigma
}^{-1}$. Note also that
\[
\bigl\|\hat{\Sigma}^{-1}-\Sigma^{-1}\bigr\|=\bigl\|\hat{
\Sigma}^{-1} (\hat{\Sigma}-\Sigma)\Sigma^{-1}\bigr\|\le\bigl\|\hat{
\Sigma}^{-1}\bigr\|\|\hat{\Sigma}-\Sigma\|\bigl\|\Sigma^{-1}\bigr\|.\vadjust{\goodbreak}
\]
Therefore,
\begin{eqnarray*}
\bE \biggl(\|\hat{\Omega}-\Omega\|^2\ind \biggl\{
\lambda_{\min}(\hat{\Sigma})\ge\frac{1}{2}\lambda_{\min}(
\Sigma) \biggr\} \biggr)&\le& 4\|\Omega\|^{2}\bE\|\hat{\Sigma }-\Sigma
\|^2
\\
&\le& C\min \biggl\{n^{-2\alpha/(2\alpha+1)}+\frac{\log
p}{n},\frac{p}{n} \biggr
\}
\end{eqnarray*}
by Theorem~\ref{thmain}. On the other hand,
\begin{eqnarray*}
&&\bE \bigl(\|\hat{\Omega}-\Omega\|^2\ind \bigl\{
\lambda_{\min}(\hat{\Sigma})< \tfrac{1}{2}\lambda_{\min}(
\Sigma) \bigr\} \bigr)
\\
&&\qquad\le \bE \bigl( \bigl(\|\hat{\Omega}\|+\|\Omega\| \bigr)^2\ind \bigl\{
\lambda_{\min}(\hat{\Sigma})< \tfrac{1}{2}\lambda_{\min}(
\Sigma) \bigr\} \bigr)
\\
&&\qquad\le\bigl(n+\|\Omega\|\bigr)^2\bP \bigl\{\lambda_{\min}(\hat{
\Sigma})< \tfrac{1}{2}\lambda_{\min}(\Sigma) \bigr\}.
\end{eqnarray*}
Note that
\[
\bP \bigl\{\lambda_{\min}(\hat{\Sigma})< \tfrac{1}{2}
\lambda_{\min
}(\Sigma) \bigr\}\le\bP \bigl\{\|\hat{\Sigma}-\Sigma\|>
\tfrac{1}{2}\lambda_{\min
}(\Sigma) \bigr\}.
\]
It suffices to show that
\[
n^2\bP \biggl\{\|\hat{\Sigma}-\Sigma\|>\frac{1}{2}
\lambda_{\min
}(\Sigma) \biggr\}\le C\min \biggl\{n^{-2\alpha/(2\alpha+1)}+
\frac{\log p}{n},\frac{p}{n} \biggr\}.
\]

Consider first the case when $p$ is large. More specifically, let
\[
p> n \bigl(48\lambda_0M^2 \bigr)^{-2}.
\]
As shown in the proof of Theorem~\ref{thmain},
\[
\bP \bigl\{\|\hat{\Sigma}-\Sigma\|>\tfrac{1}{2}\lambda_{\min
}(
\Sigma) \bigr\}\le4p^{-4}.
\]
It is not hard to see that this implies the desired claim.

Now consider the case when
\[
p\le n \bigl(48\lambda_0M^2 \bigr)^{-2}.
\]
Observe that for each $B=I\times J\in\Bcal$,
\begin{eqnarray*}
\|\hat{\Sigma}_B-\bar{\Sigma}_B\|&\le&
\lambda_0 \bigl(\|\bar{\Sigma}_{I\times I}\|\|\bar{
\Sigma}_{J\times J}\| \bigr)^{1/2} \bigl(n^{-1} \bigl(d(B)+\log p
\bigr) \bigr)^{1/2}
\\
&\le& \lambda_0\|\bar{\Sigma}\| \bigl(n^{-1} \bigl(d(B)+
\log p \bigr) \bigr)^{1/2}.
\end{eqnarray*}
It can then be deduced from the norm compression inequality, in a
similar spirit as before, that
\begin{eqnarray*}
\|\hat{\Sigma}-\bar{\Sigma}\|&\le& \sum_{l}\bigl\|S(
\hat{\Sigma}-\bar{\Sigma},l)\bigr\|
\\
&\le&12\lambda_0\|\bar{\Sigma}\|\sum_{l}
\bigl(n^{-1} \bigl(2^{l-1}k_0+\log p \bigr)
\bigr)^{1/2}
\\
&\le&12\lambda_0\|\bar{\Sigma}\|(p/n)^{1/2}.
\end{eqnarray*}
By the triangle inequality,
\[
\|\hat{\Sigma}-\Sigma\|\le\|\bar{\Sigma}-\Sigma\|+\|\hat{\Sigma }-\bar{\Sigma}
\|,\vadjust{\goodbreak}
\]
and
\[
\|\bar{\Sigma}\|\le\|\bar{\Sigma}-\Sigma\|+\|\Sigma\|.
\]
Under the event that
\[
\|\bar{\Sigma}-\Sigma\|>\frac{(1/2)\lambda_{\min}(\Sigma
)-12\lambda_0(p/n)^{1/2}\lambda_{\max}(\Sigma)}{1+12\lambda_0(p/n)^{1/2}}\ge\frac{1}{5}
\lambda_{\min}(\Sigma),
\]
we have
\[
\|\hat{\Sigma}-\Sigma\|>\tfrac{1}{2}\lambda_{\min}(\Sigma).
\]
Now by Lemma~\ref{letail},
\[
\bP \biggl\{\|\hat{\Sigma}-\Sigma\|>\frac{1}{2}\lambda_{\min
}(
\Sigma) \biggr\}\le\bP \biggl\{\|\bar{\Sigma}-\Sigma\|>\frac{1}{5}
\lambda_{\min
}(\Sigma) \biggr\}\le\exp \biggl(-\frac{cn\lambda_{\min
}^2(\Sigma)}{\lambda_{\max
}^2(\Sigma)}
\biggr),
\]
for some constant $c>0$, which concludes the proof.

\subsection{\texorpdfstring{Proof of Lemma \protect\ref{leblocknorm}}{Proof of Lemma 3}}

The proof relies on the following simple observation.

\begin{lemma}
\label{lemdist}
For any $B\in\Bcal$ with dimension $d(B)\ge4k_0$,
\[
\min_{(i,j)\in B} |i-j|\ge d(B).
\]
\end{lemma}

\begin{pf} Note that for any $B\in\Bcal$, there
exists an
integer $r>0$ such that $d(B)=2^{r-1}k_0$. We proceed by induction on
$r$. When $r=3$, it is clear by construction, blocks of size
$4k_0\times4k_0$ are at least one $2k_0\times 2k_0$ block away from the
diagonal. See Figure~\ref{figblock} also. This implies that the
statement is true for $r=3$. From~$r+1$ to $r+2$, one simply observes
that all blocks of size $2^{r+1}k_0\times2^{r+1}k_0$ is at least one
$2^rk_0\times2^rk_0$ block away from blocks of size $2^{r-1}k_0\times
2^{r-1}k_0$. Therefore,
\[
\min_{(i,j)\in B} |i-j|\ge2^rk_0+2^rk_0=2^{r+1}k_0,
\]
which implies the desired statement.
\end{pf}

We are now in position to prove Lemma~\ref{leblocknorm} which states
that big blocks of the covariance matrix are small in spectral norm.
Recall that the matrix $\ell_1$ norm is defined as
\[
\|A\|_{\ell_1}=\sup_{\bx\in\bR^p:\|\bx\|_{\ell_1}=1}\|A\bx\|_{\ell
_1}=
\max_{1\le j\le n}\sum_{i=m}^p
|a_{ij}|,
\]
for an $m\times n$ matrix $A=(a_{ij}))_{1\le i\le m, 1\le j\le n}$.
Similarly the matrix $\ell_\infty$ norm is defined as
\[
\|A\|_{\ell_\infty}=\sup_{\bx\in\bR^p:\|\bx\|_{\ell_\infty
}=1}\|A\bx\|_{\ell_\infty}=
\max_{1\le i\le m}\sum_{j=1}^n
|a_{ij}|.\vadjust{\goodbreak}
\]
It is well known [see, e.g., \citet{gv96}] that
\[
\|A\|^2\le\|A\|_{\ell_1}\|A\|_{\ell_\infty}.
\]
Immediately from Lemma~\ref{lemdist}, we have
\[
\|\Sigma_B\|_{\ell_1}, \|\Sigma_B
\|_{\ell_\infty}\le\max_{1\le
i\le
p}\sum_{j: |j-i|\ge2^rk_0} |
\sigma_{ij}|\le Md(B)^{-\alpha},
\]
which implies
$
\|\Sigma_B\|\le Md(B)^{-\alpha}$.

\subsection{\texorpdfstring{Proof of Lemma \protect\ref{letune}}{Proof of Lemma 4}}
For any $I\in\Ical$, write $Z_I=\Sigma_{I\times I}^{-1/2}Y$. Then the
entries of $Z_I$ are independent standard normal random variables. From
the concentration bounds on the random matrices [see, e.g., \citet
{ds01}], we have
\[
1-\frac{\sqrt{\card(I)}+t}{\sqrt{n}}\le\lambda_{\min}^{1/2}(\bar {\Sigma
}_{Z_I})\le\lambda_{\max}^{1/2}(\bar{
\Sigma}_{Z_I})\le1+\frac{\sqrt{\card
(I)}+t}{\sqrt{n}}
\]
with probability at least $1-2\exp(-t^2/2)$ where $\bar{\Sigma}_{Z_I}$
is the sample covariance matrix of $Z_I$. Applying the union bound to
all $I\in\Ical$ yields that with probability at least $1-2p^2\exp
(-t^2/2)$, for all $I$
\[
1-\frac{\sqrt{\card(I)}+t}{\sqrt{n}}\le\lambda_{\min}^{1/2}(\bar {\Sigma
}_{Z_I})\le\lambda_{\max}^{1/2}(\bar{
\Sigma}_{Z_I})\le1+\frac{\sqrt{\card
(I)}+t}{\sqrt{n}}.
\]
%

Observe that
$
\bar{\Sigma}_{I\times I}=\Sigma_{I\times I}^{1/2}\bar{\Sigma
}_{Z_I}\Sigma_{I\times I}^{1/2}.
$
Thus
\[
\lambda_{\min}(\bar{\Sigma}_{Z_I})\lambda_{\max}(
\Sigma_{I\times
I})\le\lambda_{\max}(\bar{\Sigma}_{I\times I})\le
\lambda_{\max}(\bar{\Sigma}_{Z_I})\lambda_{\max}(
\Sigma_{I\times I}),
\]
which implies the desired statement.

\printaddresses


\begin{thebibliography}{27}

\bibitem[\protect\citeauthoryear{Banerjee, El Ghaoui and d'Aspremont}{2008}]{bgd08}
%
\begin{barticle}[mr]
\bauthor{\bsnm{Banerjee},~\bfnm{Onureena}\binits{O.}},
\bauthor{\bsnm{El~Ghaoui},~\bfnm{Laurent}\binits{L.}} \AND
\bauthor{\bsnm{d'Aspremont},~\bfnm{Alexandre}\binits{A.}}
(\byear{2008}).
\btitle{Model selection through sparse maximum likelihood estimation for
multivariate {G}aussian or binary data}.
\bjournal{J.~Mach. Learn. Res.}
\bvolume{9}
\bpages{485--516}.
\bid{issn={1532-4435}, mr={2417243}}
\bptok{imsref}%
\end{barticle}
%
\endbibitem

\bibitem[\protect\citeauthoryear{Bickel and Levina}{2008a}]{bl08a}
%
\begin{barticle}[mr]
\bauthor{\bsnm{Bickel},~\bfnm{Peter~J.}\binits{P.~J.}} \AND
\bauthor{\bsnm{Levina},~\bfnm{Elizaveta}\binits{E.}}
(\byear{2008}a).
\btitle{Regularized estimation of large covariance matrices}.
\bjournal{Ann. Statist.}
\bvolume{36}
\bpages{199--227}.
\bid{doi={10.1214/009053607000000758}, issn={0090-5364}, mr={2387969}}
\bptok{imsref}%
\end{barticle}
%
\endbibitem

\bibitem[\protect\citeauthoryear{Bickel and Levina}{2008b}]{bl08b}
%
\begin{barticle}[mr]
\bauthor{\bsnm{Bickel},~\bfnm{Peter~J.}\binits{P.~J.}} \AND
\bauthor{\bsnm{Levina},~\bfnm{Elizaveta}\binits{E.}}
(\byear{2008}b).
\btitle{Covariance regularization by thresholding}.
\bjournal{Ann. Statist.}
\bvolume{36}
\bpages{2577--2604}.
\bid{doi={10.1214/08-AOS600}, issn={0090-5364}, mr={2485008}}
\bptok{imsref}%
\end{barticle}
%
\endbibitem

\bibitem[\protect\citeauthoryear{B{\"o}r{\"o}czky and Wintsche}{2005}]{bw05}
%
\begin{bmisc}[auto:STB|2012/09/06|07:18:59]
\bauthor{\bsnm{B{\"o}r{\"o}czky},~\bfnm{K.}\binits{K.}} \AND
\bauthor{\bsnm{Wintsche},~\bfnm{G.}\binits{G.}}
(\byear{2005}).
\bhowpublished{Covering the sphere by equal spherical balls. Available at
\url{http://www.renyi.hu/\textasciitilde carlos/spherecover.ps}.}
\bptok{imsref}%
\end{bmisc}
%
\endbibitem

\bibitem[\protect\citeauthoryear{Cai}{1999}]{c99}
%
\begin{barticle}[mr]
\bauthor{\bsnm{Cai},~\bfnm{T.~Tony}\binits{T.~T.}}
(\byear{1999}).
\btitle{Adaptive wavelet estimation: A block thresholding and oracle inequality
approach}.
\bjournal{Ann. Statist.}
\bvolume{27}
\bpages{898--924}.
\bid{doi={10.1214/aos/1018031262}, issn={0090-5364}, mr={1724035}}
\bptok{imsref}%
\end{barticle}
%
\endbibitem

\bibitem[\protect\citeauthoryear{Cai and Liu}{2011}]{cl11}
%
\begin{barticle}[mr]
\bauthor{\bsnm{Cai},~\bfnm{Tony}\binits{T.}} \AND
\bauthor{\bsnm{Liu},~\bfnm{Weidong}\binits{W.}}
(\byear{2011}).
\btitle{Adaptive thresholding for sparse covariance matrix estimation}.
\bjournal{J. Amer. Statist. Assoc.}
\bvolume{106}
\bpages{672--684}.
\bid{doi={10.1198/jasa.2011.tm10560}, issn={0162-1459}, mr={2847949}}
\bptok{imsref}%
\end{barticle}
%
\endbibitem

\bibitem[\protect\citeauthoryear{Cai, Liu and Luo}{2011}]{cll11}
%
\begin{barticle}[mr]
\bauthor{\bsnm{Cai},~\bfnm{Tony}\binits{T.}},
\bauthor{\bsnm{Liu},~\bfnm{Weidong}\binits{W.}} \AND
\bauthor{\bsnm{Luo},~\bfnm{Xi}\binits{X.}}
(\byear{2011}).
\btitle{A constrained {$\ell\sb1$} minimization approach to sparse precision
matrix estimation}.
\bjournal{J. Amer. Statist. Assoc.}
\bvolume{106}
\bpages{594--607}.
\bid{doi={10.1198/jasa.2011.tm10155}, issn={0162-1459}, mr={2847973}}
\bptok{imsref}%
\end{barticle}
%
\endbibitem

\bibitem[\protect\citeauthoryear{Cai, Liu and Zhou}{2011}]{clz11}
%
\begin{bmisc}[auto:STB|2012/09/06|07:18:59]
\bauthor{\bsnm{Cai},~\bfnm{T.~T.}\binits{T.~T.}},
\bauthor{\bsnm{Liu},~\bfnm{W.}\binits{W.}} \AND
\bauthor{\bsnm{Zhou},~\bfnm{H.~H.}\binits{H.~H.}}
(\byear{2011}).
\bhowpublished{Optimal estimation of large sparse precision matrices.
Unpublished manuscript}.
\bptok{imsref}%
\end{bmisc}
%
\endbibitem

\bibitem[\protect\citeauthoryear{Cai, Zhang and Zhou}{2010}]{czz10}
%
\begin{barticle}[mr]
\bauthor{\bsnm{Cai},~\bfnm{T.~Tony}\binits{T.~T.}},
\bauthor{\bsnm{Zhang},~\bfnm{Cun-Hui}\binits{C.-H.}} \AND
\bauthor{\bsnm{Zhou},~\bfnm{Harrison~H.}\binits{H.~H.}}
(\byear{2010}).
\btitle{Optimal rates of convergence for covariance matrix estimation}.
\bjournal{Ann. Statist.}
\bvolume{38}
\bpages{2118--2144}.
\bid{doi={10.1214/09-AOS752}, issn={0090-5364}, mr={2676885}}
\bptok{imsref}%
\end{barticle}
%
\endbibitem

\bibitem[\protect\citeauthoryear{Cai and Zhou}{2011}]{cz11}
%
\begin{bmisc}[auto:STB|2012/09/06|07:18:59]
\bauthor{\bsnm{Cai},~\bfnm{T.~T.}\binits{T.~T.}} \AND
\bauthor{\bsnm{Zhou},~\bfnm{H.}\binits{H.}}
(\byear{2011}).
\bhowpublished{Optimal rates of convergence for sparse covariance
matrix estimation. Technical report}.
\bptok{imsref}%
\end{bmisc}
%
\endbibitem

\bibitem[\protect\citeauthoryear{Davidson and Szarek}{2001}]{ds01}
%
\begin{bincollection}[mr]
\bauthor{\bsnm{Davidson},~\bfnm{Kenneth~R.}\binits{K.~R.}} \AND
\bauthor{\bsnm{Szarek},~\bfnm{Stanislaw~J.}\binits{S.~J.}}
(\byear{2001}).
\btitle{Local operator theory, random matrices and {B}anach spaces}.
In \bbooktitle{Handbook of the Geometry of {B}anach Spaces, {V}ol. {I}}
\bpages{317--366}.
\bpublisher{North-Holland}, \blocation{Amsterdam}.
\bid{doi={10.1016/S1874-5849(01)80010-3}, mr={1863696}}
\bptok{imsref}%
\end{bincollection}
%
\endbibitem

\bibitem[\protect\citeauthoryear{Efromovich}{1985}]{e85}
%
\begin{barticle}[auto:STB|2012/09/06|07:18:59]
\bauthor{\bsnm{Efromovich},~\bfnm{S.~Y.}\binits{S.~Y.}}
(\byear{1985}).
\btitle{Nonparametric estimation of a density of unknown smoothness}.
\bjournal{Theory Probab. Appl.}
\bvolume{30}
\bpages{557--661}.
\bptok{imsref}%
\end{barticle}
%
\endbibitem

\bibitem[\protect\citeauthoryear{El Karoui}{2008}]{e08}
%
\begin{barticle}[mr]
\bauthor{\bsnm{El~Karoui},~\bfnm{Noureddine}\binits{N.}}
(\byear{2008}).
\btitle{Operator norm consistent estimation of large-dimensional sparse
covariance matrices}.
\bjournal{Ann. Statist.}
\bvolume{36}
\bpages{2717--2756}.
\bid{doi={10.1214/07-AOS559}, issn={0090-5364}, mr={2485011}}
\bptok{imsref}%
\end{barticle}
%
\endbibitem

\bibitem[\protect\citeauthoryear{Fan, Fan and Lv}{2008}]{ffl08}
%
\begin{barticle}[mr]
\bauthor{\bsnm{Fan},~\bfnm{Jianqing}\binits{J.}},
\bauthor{\bsnm{Fan},~\bfnm{Yingying}\binits{Y.}} \AND
\bauthor{\bsnm{Lv},~\bfnm{Jinchi}\binits{J.}}
(\byear{2008}).
\btitle{High dimensional covariance matrix estimation using a factor model}.
\bjournal{J. Econometrics}
\bvolume{147}
\bpages{186--197}.
\bid{doi={10.1016/j.jeconom.2008.09.017}, issn={0304-4076}, mr={2472991}}
\bptok{imsref}%
\end{barticle}
%
\endbibitem

\bibitem[\protect\citeauthoryear{Friedman, Hastie and
Tibshirani}{2008}]{fht08}
%
\begin{barticle}[auto:STB|2012/09/06|07:18:59]
\bauthor{\bsnm{Friedman},~\bfnm{J.}\binits{J.}},
\bauthor{\bsnm{Hastie},~\bfnm{T.}\binits{T.}} \AND
\bauthor{\bsnm{Tibshirani},~\bfnm{T.}\binits{T.}}
(\byear{2008}).
\btitle{Sparse inverse covariance estimation with the graphical lasso}.
\bjournal{Biostatistics}
\bvolume{9}
\bpages{432--441}.
\bptok{imsref}%
\end{barticle}
%
\endbibitem

\bibitem[\protect\citeauthoryear{Golub and Van Loan}{1996}]{gv96}
%
\begin{bbook}[mr]
\bauthor{\bsnm{Golub},~\bfnm{Gene~H.}\binits{G.~H.}} \AND
\bauthor{\bsnm{Van~Loan},~\bfnm{Charles~F.}\binits{C.~F.}}
(\byear{1996}).
\btitle{Matrix Computations},
\bedition{3rd} ed.
\bpublisher{Johns Hopkins Univ. Press}, \blocation{Baltimore, MD}.
\bid{mr={1417720}}
\bptok{imsref}%
\end{bbook}
%
\endbibitem

\bibitem[\protect\citeauthoryear{Huang et al.}{2006}]{hlpl06}
%
\begin{barticle}[mr]
\bauthor{\bsnm{Huang},~\bfnm{Jianhua~Z.}\binits{J.~Z.}},
\bauthor{\bsnm{Liu},~\bfnm{Naiping}\binits{N.}},
\bauthor{\bsnm{Pourahmadi},~\bfnm{Mohsen}\binits{M.}} \AND
\bauthor{\bsnm{Liu},~\bfnm{Linxu}\binits{L.}}
(\byear{2006}).
\btitle{Covariance matrix selection and estimation via penalised normal
likelihood}.
\bjournal{Biometrika}
\bvolume{93}
\bpages{85--98}.
\bid{doi={10.1093/biomet/93.1.85}, issn={0006-3444}, mr={2277742}}
\bptok{imsref}%
\end{barticle}
%
\endbibitem

\bibitem[\protect\citeauthoryear{Lam and Fan}{2009}]{lf09}
%
\begin{barticle}[mr]
\bauthor{\bsnm{Lam},~\bfnm{Clifford}\binits{C.}} \AND
\bauthor{\bsnm{Fan},~\bfnm{Jianqing}\binits{J.}}
(\byear{2009}).
\btitle{Sparsistency and rates of convergence in large covariance matrix
estimation}.
\bjournal{Ann. Statist.}
\bvolume{37}
\bpages{4254--4278}.
\bid{doi={10.1214/09-AOS720}, issn={0090-5364}, mr={2572459}}
\bptok{imsref}%
\end{barticle}
%
\endbibitem

\bibitem[\protect\citeauthoryear{Laurent and Massart}{2000}]{lm00}
%
\begin{barticle}[mr]
\bauthor{\bsnm{Laurent},~\bfnm{B.}\binits{B.}} \AND
\bauthor{\bsnm{Massart},~\bfnm{P.}\binits{P.}}
(\byear{2000}).
\btitle{Adaptive estimation of a quadratic functional by model selection}.
\bjournal{Ann. Statist.}
\bvolume{28}
\bpages{1302--1338}.
\bid{doi={10.1214/aos/1015957395}, issn={0090-5364}, mr={1805785}}
\bptok{imsref}%
\end{barticle}
%
\endbibitem

\bibitem[\protect\citeauthoryear{Ledoit and Wolf}{2004}]{lw04}
%
\begin{barticle}[mr]
\bauthor{\bsnm{Ledoit},~\bfnm{Olivier}\binits{O.}} \AND
\bauthor{\bsnm{Wolf},~\bfnm{Michael}\binits{M.}}
(\byear{2004}).
\btitle{A well-conditioned estimator for large-dimensional covariance
matrices}.
\bjournal{J. Multivariate Anal.}
\bvolume{88}
\bpages{365--411}.
\bid{doi={10.1016/S0047-259X(03)00096-4}, issn={0047-259X}, mr={2026339}}
\bptok{imsref}%
\end{barticle}
%
\endbibitem

\bibitem[\protect\citeauthoryear{Ravikumar et al.}{2011}]{rwry}
%
\begin{barticle}[mr]
\bauthor{\bsnm{Ravikumar},~\bfnm{Pradeep}\binits{P.}},
\bauthor{\bsnm{Wainwright},~\bfnm{Martin~J.}\binits{M.~J.}},
\bauthor{\bsnm{Raskutti},~\bfnm{Garvesh}\binits{G.}} \AND
\bauthor{\bsnm{Yu},~\bfnm{Bin}\binits{B.}}
(\byear{2011}).
\btitle{High-dimensional covariance estimation by minimizing {$\ell
\sb
1$}-penalized log-determinant divergence}.
\bjournal{Electron. J. Stat.}
\bvolume{5}
\bpages{935--980}.
\bid{doi={10.1214/11-EJS631}, issn={1935-7524}, mr={2836766}}
\bptok{imsref}%
\end{barticle}
%
\endbibitem

\bibitem[\protect\citeauthoryear{Rocha, Zhao and Yu}{2008}]{rzy08}
%
\begin{bmisc}[auto:STB|2012/09/06|07:18:59]
\bauthor{\bsnm{Rocha},~\bfnm{G.}\binits{G.}},
\bauthor{\bsnm{Zhao},~\bfnm{P.}\binits{P.}} \AND
\bauthor{\bsnm{Yu},~\bfnm{B.}\binits{B.}}
(\byear{2008}).
\bhowpublished{A path following algorithm for sparse pseudo-likelihood
inverse covariance estimation. Technical report, Dept. Statistics,
Univ. California, Berkeley}.
\bptok{imsref}%
\end{bmisc}
%
\endbibitem

\bibitem[\protect\citeauthoryear{Rothman, Levina and Zhu}{2009}]{rlz09}
%
\begin{barticle}[mr]
\bauthor{\bsnm{Rothman},~\bfnm{Adam~J.}\binits{A.~J.}},
\bauthor{\bsnm{Levina},~\bfnm{Elizaveta}\binits{E.}} \AND
\bauthor{\bsnm{Zhu},~\bfnm{Ji}\binits{J.}}
(\byear{2009}).
\btitle{Generalized thresholding of large covariance matrices}.
\bjournal{J. Amer. Statist. Assoc.}
\bvolume{104}
\bpages{177--186}.
\bid{doi={10.1198/jasa.2009.0101}, issn={0162-1459}, mr={2504372}}
\bptok{imsref}%
\end{barticle}
%
\endbibitem

\bibitem[\protect\citeauthoryear{Rothman et al.}{2008}]{rblz08}
%
\begin{barticle}[mr]
\bauthor{\bsnm{Rothman},~\bfnm{Adam~J.}\binits{A.~J.}},
\bauthor{\bsnm{Bickel},~\bfnm{Peter~J.}\binits{P.~J.}},
\bauthor{\bsnm{Levina},~\bfnm{Elizaveta}\binits{E.}} \AND
\bauthor{\bsnm{Zhu},~\bfnm{Ji}\binits{J.}}
(\byear{2008}).
\btitle{Sparse permutation invariant covariance estimation}.
\bjournal{Electron. J. Stat.}
\bvolume{2}
\bpages{494--515}.
\bid{doi={10.1214/08-EJS176}, issn={1935-7524}, mr={2417391}}
\bptok{imsref}%
\end{barticle}
%
\endbibitem

\bibitem[\protect\citeauthoryear{Yuan}{2010}]{y10}
%
\begin{barticle}[mr]
\bauthor{\bsnm{Yuan},~\bfnm{Ming}\binits{M.}}
(\byear{2010}).
\btitle{High dimensional inverse covariance matrix estimation via linear
programming}.
\bjournal{J. Mach. Learn. Res.}
\bvolume{11}
\bpages{2261--2286}.
\bid{issn={1532-4435}, mr={2719856}}
\bptok{imsref}%
\end{barticle}
%
\endbibitem

\bibitem[\protect\citeauthoryear{Yuan and Lin}{2007}]{yl07}
%
\begin{barticle}[mr]
\bauthor{\bsnm{Yuan},~\bfnm{Ming}\binits{M.}} \AND
\bauthor{\bsnm{Lin},~\bfnm{Yi}\binits{Y.}}
(\byear{2007}).
\btitle{Model selection and estimation in the {G}aussian graphical model}.
\bjournal{Biometrika}
\bvolume{94}
\bpages{19--35}.
\bid{doi={10.1093/biomet/asm018}, issn={0006-3444}, mr={2367824}}
\bptok{imsref}%
\end{barticle}
%
\endbibitem

\end{thebibliography}
\end{document}